\documentstyle[12pt]{article}

\newtheorem{Thm}{Theorem}[section]

\newtheorem{Cor}[Thm]{Corollary}
\newtheorem{Lem}[Thm]{Lemma}
\newtheorem{Prop}[Thm]{Proposition}
\newtheorem{Def}[Thm]{Definition}
\newtheorem{Rmk}[Thm]{Remark}
\newtheorem{Obs}[Thm]{Observation}

\newcommand{\im}{\mbox{$\Rightarrow$}}
\newcommand{\eq}{\mbox{$\Leftrightarrow$}}

\newcommand{\supp}{{\rm supp}\,}
\newcommand{\subl}{{_{\mbox {\scriptsize $\ell$}}}}
\newcommand{\subk}{{_{\mbox {\scriptsize $k$}}}}
\renewcommand{\span}{{\rm span}\,}
\newcommand{\ran}{{\rm ran}\,}

\def \norm{| \! | \! |}
\def \N{{\rm I\kern-2.9pt{\rm N}}}
\def \R{{\rm I\kern-2.9pt{\rm R}}}
\def \K{{\rm I\kern-2.9pt{\rm K}}}

\begin{document}

\title{Distorting Mixed Tsirelson Spaces }

\author{G. Androulakis \thanks{Part of this paper also appears in the first
author's Ph.D. thesis which is being prepared  under the supervision of 
Prof. H. Rosenthal at the University of Texas at Austin.}
\and E. Odell \thanks{Research supported by NSF and TARP}}

\date{}
\maketitle

\noindent
{\bf Abstract:}
Any regular mixed Tsirelson space $T(\theta_n,S_n)_{\N}$ for which $\frac{
\theta_n}{\theta^n} \rightarrow 0$, where $\theta=\lim_n \theta_n^{1/n}$, is
shown to be arbitrarily distortable. Certain asymptotic $\ell_1$ constants for
those and other mixed Tsirelson spaces are calculated. Also a combinatorial
result on the Schreier families $(S_{\alpha})_{\alpha < \omega_1}$ is proved
and an application is given to show that for every Banach space $X$ with a 
basis $(e_i)$, the two $\Delta$-spectrums $\Delta(X)$ and $\Delta(X,(e_i))$
coincide.

\noindent
\section{Introduction} \label{S:intro}
A Banach space $X$ with basis $(e_i)$ is asymptotic $\ell_1$ if there exists
$\delta >0$ such that
for all $n$ and block bases $(x_i)_1^n$ of $(e_i)_n^{\infty}$,
\begin{equation} \label{E:asl1}
\| \sum_{i=1}^n x_i \| \geq \delta \sum_{i=1}^n\| x_i \|.
\end{equation}
Such a space need not contain $\ell_1$ as witnessed by Tsirelson's famous 
space $T$. The complexity of the asymptotic $\ell_1$ structure within $X$ 
can be
measured by certain constants $\delta_{\alpha}(e_i)$ for $\alpha<\omega_1$.
$\delta_1(e_i)$ is the largest $\delta >0$ satisfying (\ref{E:asl1}) above. 
Subsequent $\delta_{\alpha}$'s are defined by a similar formula where $(x_i)_1
^n$ ranges over ``$\alpha$-admissible'' block bases (all terms are precisely 
defined in section~\ref{S:prel}). These notions were developed in \cite{OTW}
where, in addition, $\delta_{\alpha}(y_i)$ was considered, for a block basis 
$(y_i)$ of $(e_i)$. In this setting, $(y_i)$ becomes the reference frame and
one naturally has $\delta_{\alpha}(y_i) \geq \delta_{\alpha}(e_i)$. These 
constants can perhaps increase by passing to further block bases and this 
leads to the notion of the $\Delta$-spectrum of $X$, $\Delta(X)$. Roughly, 
$\Delta(X)$ is the set of all $\gamma =(\gamma_{\alpha})_{\alpha< \omega_1}$
where $\gamma_{\alpha}$ is the stabilization of $\delta
_{\alpha}(y_i)$ for $(y_i)$ some block basis of $(e_i)$. Alternatively 
by  keeping $(e_i)$ as the reference frame,
in a similar manner we obtain $\Delta(X,(e_i))$. In section~\ref{S:Delta} we
prove that these two notions coincide, $\Delta(X)=\Delta(X,(e_i))$.

\bigskip
\noindent
Argyros and Deliyanni \cite{AD} constructed the first example of an asymptotic
$\ell_1$ arbitrarily distortable Banach space by
constructing ``mixed Tsirelson spaces'' and
proving that such spaces can be arbitrarily distortable. In 
section~\ref{S:mixedtsirelson} we consider the simplest class of mixed
Tsirelson spaces
$X= T(\theta_n,S_n)_{n \in \N}$ where $\theta_n \rightarrow 0$ and $\sup_n 
\theta_n <1$. These are reflexive asymptotic $\ell_1$ spaces having a 
1-unconditional basis $(e_i)$. Also we may assume $\theta \equiv 
\theta_n^{1/n} $ exists. We prove that if $\frac{\theta_n}{\theta^n} 
\rightarrow 0$ then $X$ is arbitrarily distortable. In particular, this happens
if $\theta=1$. Thus, for example, $T(\frac{1}{n+1},S_n)_{\N}$ is an arbitrarily
distortable space. We also calculate the asymptotic constants $\ddot{\delta}
_{\alpha}(X)$ for these spaces along with the spectral index $I_{\Delta}(X)$.
$\ddot{\delta}_{\alpha}(X)$ is the supremum of $\delta_{\alpha}((x_i),\mid 
\cdot \mid)$ under all equivalent norms on $X$ and $I_{\Delta}(X)$ is the first
ordinal $\alpha$ for which $\ddot{\delta}_{\alpha}(X)<1$.

\noindent
\section{Preliminaries} \label{S:prel}
$X, Y, Z, \ldots$ shall denote separable infinite dimensional Banach spaces.
{\em All the spaces we consider will have bases.} 
\/Every Banach space with a basis can be viewed as the completion of $c_{00}$
(the linear space of finitely supported real valued sequences) 
under a certain norm. $(e_i)$
will denote the unit vector basis for $c_{00}$ and whenever a Banach space 
$(X, \| \cdot \|)$ with
a basis is regarded as the completion of $(c_{00}, \| \cdot \|)$,
  $(e_i)$ will denote this (normalized)  basis.
If $x \in c_{00}$ and $E \subseteq \N$, $Ex \in c_{00}$ is the restriction 
of $x$ to $E$;
$Ex(j)=x(j)$ if $j \in E$ and 0 otherwise.
Also the {\em support} \/ of $x$, $\supp(x)$, (w.r.t. $(e_i)$) is the set 
$\{ j \in \N : x(j) \not = 0\}$. The {\em range}\/ of $x$, $\ran(x)$, 
(w.r.t. $(e_i)$) is the smallest
interval which contains $\supp(x)$. If $x$, $x_1$, $x_2, \ldots $ are 
vectors in $X$ (and $k \in \N$) then we say that $x$ is an average of $(x_i)_i$
(of length $k$) if there exists $F \subset \N$ with $x=\frac{1}{\mid F \mid} 
\sum_{i \in F} x_i $ (and $\mid F \mid =k$).  
We say that a sequence $(y_i)$ is a 
{\em (convex) block sequence of $(x_i)$} \/ if for all $i$, $y_i= \sum_{j=m_i}
^{m_{i+1}-1}\alpha_j x_j$ for some sequence  $m_1<m_2< \ldots$, 
of integers and $(\alpha_j)_{j \in \N} \subset \R$ (resp. \ with $\alpha_j \geq 
0$ for all $j$, and $\sum_{j=m_i}^{m_{i+1}-1}\alpha_j=1$ for all $i$).
If $(x_i)$ is a block basis of
$(y_i)$ we write $(x_i) \prec (y_i)$. $X \prec Y$ shall mean that $X$ has a 
basis which is a block basis of a certain basis for $Y$, when the given 
bases are understood. For $\lambda  >1$, $(X, \| \cdot \|)$ is 
{\em $\lambda$-distortable} \/ if there exists an equivalent norm $\mid \cdot 
\mid$ on $X$ so that for all $Y \prec X$ 
\[ D(Y, \mid \cdot \mid) \equiv \sup \{ \frac{\mid y \mid}{\mid z \mid} : 
y, z \in Y, \| y \| = \| z \| =1 \} \geq \lambda. \] 
$X$ is {\em distortable} \/ if it is $\lambda$-distortable for some
$\lambda >1$ and {\em arbitrarily distortable} \/ if it is $\lambda$-distortable
for all $\lambda >1$. $X$ is of {\em $D$-bounded distortion} \/ if for all
equivalent norms $\mid \cdot \mid$ on $X$ and for all $Z \prec X$ there exists
$Y \prec Z$ with $D(Y, \mid \cdot \mid) \leq D$. Note that if $C= \inf \{
\frac{\| y \|}{\mid y \mid}: y \in Y, y \not =0 \}$ then
\begin{equation} \label{E:equiv}
C \mid y \mid \leq \| y \| \leq D(Y, \mid \cdot \mid) C \mid y \mid 
\mbox{for all } y \in Y.
\end{equation}
For more information on distortion we 
recommend the reader consult the following papers: \cite{S}, \cite{MT},
\cite{OS1}, \cite{OS2}, \cite{OS3}, \cite{Ma}, \cite{T}, \cite{OTW}.  

\bigskip
\noindent
Asymptotic $\ell_1$ Banach spaces are defined by  (\ref{E:asl1}) in 
section~\ref{S:intro} (for another approach to asymptotic structure see
\cite{MMT}). 
These spaces were studied in \cite{OTW} where certain asymptotic
constants were introduced. We shall recall the relevant definitions but first
we need to recall the definition of the Schreier sets $S_{\alpha}$, $\alpha
< \omega _1$ \cite{AA}. For $F, G \subset \N$,  we write $F < G$ when
$\max (F) < \min (G)$ or one of them is empty, and we write $n \leq F$ instead
of $\{ n \} \leq F$. Also for $x, y \in c_{00}$, $x<y$ means $\ran(x)< \ran(y)$.

\begin{Def} $S_0 = \{ \{ n \} : n \in \N \} \cup \{ \emptyset \}$. 
If $\alpha < \omega _1$ and $S_{
\alpha}$ has been defined, \\ 
$S_{\alpha +1} = \{ \cup_1 ^n F_i : n \in \N , 
n \leq F_1 <F_2< \cdots <F_n \mbox{ and } F_i \in S_{ \alpha } \mbox{ for } 
1 \leq i \leq n \}. $\\ 
If $\alpha$ is a limit ordinal choose $\alpha _n 
\nearrow \alpha$ and set $S_{\alpha} = \{ F : n \leq F \in S_{\alpha _n} 
\mbox{ for some } n \}$. 
\end{Def}

\noindent
If $(E_i)_1^\ell$ is a finite sequence of non-empty subsets of $\N$ and 
$\alpha < \omega_1$ then we say that {\em $(E_i)_1^\ell$ is 
$\alpha$-admissible}\/ if $E_1< \cdots < E_\subl$ and $(\min \ E_i)_1^\ell
\in S_{\alpha}$. \/
If $(e_i)$ is a basic sequence and $(x_i)_1 ^\ell \prec
(e_i)$ then $(x_i)_1 ^\ell$ is {\em $\alpha$-admissible with respect to 
$(e_i)$}
\/ if $( \ran(x_i) )_1 ^\ell$ is $\alpha$-admissible where the range of $x$, 
$\ran(x)$, is w.r.t $(e_i)$. \/ If $x \in \span (x_i)$, 
then $x$ is {\em $\alpha$-admissible w.r.t. $(x_i)$} \/ if $\supp(x)$ (w.r.t. 
$(x_i)$) $\in S_\alpha$. Also if $x \in \span(x_i)$ then {\em  $x$
is a 1-admissible average of $(x_i)$ w.r.t. $(e_i)$} \/ if there exists a 
finite
set $F \subset \N$ such that $x= \frac{1}{ \mid F \mid} \sum_{i \in F} x_i$
and $(x_i)_{i \in F}$ is 1-admissible w.r.t. $(e_i)$. Note that if $x$
is a 1-admissible average of $(x_i)$ w.r.t. $(e_i)$ and for some $\alpha < 
\omega_1$ each $x_i$ is $\alpha$-admissible w.r.t. $(e_i)$ then $x$ is $\alpha
+1$-admissible w.r.t. $(e_i)$.
Thus if $(x_i)$ is a basis for $X$ then $X$ is asymptotic $\ell_1$ iff 
\begin{eqnarray*}
0< \delta _1 (x_i) \equiv \delta _1(X) \equiv \delta_1(X, \| \cdot \|)= 
\sup \! \! \! & \{ & \! \! \! \delta 
\geq 0: \| \sum _1^n y_i \| \geq \delta \sum _1 ^n \| y_i \| \mbox{ whenever } 
(y_i)_1^n \prec(x_i) \\
& & \! \! \! \mbox{and } (y_i)_1 ^n \mbox{ is 1-admissible w.r.t. } (x_i) \}.
\end{eqnarray*}
In \cite{OTW} this definition was extended as follows: 
For $\alpha < \omega_1$
\begin{eqnarray*}
\delta _{\alpha} (x_i) \equiv \delta _{\alpha}(X) \equiv \delta_{\alpha}
(X, \| \cdot \|)= \sup \! \! \! & \{ & \! \! 
\! \delta \geq 0: \| \sum_1 ^n y_i \| \geq \delta \sum _1 ^n \| y_i \| 
\mbox{ whenever }(y_i)_1^n \prec (x_i) \\
& & \! \! \!  \mbox{and } (y_i)_1 ^n \mbox{ is } \alpha \mbox{-admissible  
w.r.t. } (x_i) \}.
\end{eqnarray*}

\begin{Obs} \label{O:equiv}
Note that if we have two equivalent norms $\| \cdot \|$, $\norm \cdot \norm$ on
$X$ and for some $c, C>0$, $c \norm x \norm \leq \| x \| \leq C \norm x \norm$
for all $x \in X$, then for all $\alpha < \omega_1$,
\[ \frac{c}{C} \delta_{\alpha}(X, \norm \cdot \norm) \leq \delta_{\alpha}(X,
\| \cdot \|) \leq \frac{C}{c} \delta_{\alpha}(X, \norm \cdot \norm). \]
\end{Obs}
In problems of distortion one is 
concerned with block bases and equivalent norms. Thus we also consider 
\cite{OTW} 
\[ \begin{array}{ll}
\dot{\delta}_{\alpha}(x_i) &  = \sup \{ \delta _{\alpha} (y_i) :
(y_i) \prec (x_i) \} \mbox{ and } \\ 
\ddot{\delta}_{\alpha}(x_i) &  = \sup \{ \dot{ \delta}
_{\alpha}((x_i), \mid \cdot \mid ) : \mid \cdot \mid \mbox{ is an equivalent
norm on } X \}. 
\end{array} \]
If $(y_i) \prec (x_i)$ then $\delta_{\alpha} (y_i) \geq \delta _{\alpha}(x_i)$.
This is because each $S_{\alpha}$ is {\em spreading} \/ (if $(n_i)_1 ^k 
\in S_{\alpha}$
and $m_1< \cdots <m_\subk$ with $n_i \leq m_i$ for all $i=1, \ldots , k$, then
$(m_i)_1 ^k \in S_{\alpha}$). This leads to the following definition 
\cite{OTW}.

\begin{Def} A basic sequence $(y_i)$ $\Delta$-stabilizes  $\gamma = (\gamma
_{\alpha})_{\alpha < \omega _1} \subseteq \R$ if there exists 
$\varepsilon _n \searrow 0$
so that for all $\alpha < \omega _1$ there exists $m \in \N$ so that for all
$n \geq m$ if $(z_i) \prec (y_i)_n ^{\infty}$ then $\mid \delta_{\alpha}(z_i)-
\gamma_{\alpha} \mid < \varepsilon _n$. 
\end{Def}

\noindent
{\bf Remark} It is automatic from the definition that if $(y_i)$ 
$\Delta$-stabilizes $\gamma$ then for all $\alpha < \omega_1$, 
$\gamma_{\alpha}= \sup \{ \delta_{\alpha}(z_i): (z_i) \prec (y_i) \}$. 
Furthermore if $(z_i) \prec (y_i)$ then $(z_i)$ $\Delta$-stabilizes $\gamma$.

\bigskip
\noindent
It is shown in \cite{OTW} that if $X$ has a basis $(x_i)$ and $(y_i) \prec
(x_i)$ then there exists $(z_i) \prec (y_i)$ and $\gamma = (\gamma _{\alpha})
_{\alpha < \omega _1}$ so that $(z_i)$ $ \Delta$-stabilizes $\gamma$.   

\begin{Def}
Let $X$ have a basis $(x_i)$. The $\Delta$-spectrum of $X$, $\Delta (X)$, is 
defined to be the set of all $\gamma$'s so that $(y_i)$ stabilizes  
$\gamma$ for some $(y_i) \prec (x_i)$. 
We also define $\ddot{\Delta} (X) = \cup \{ \Delta(X, \mid \cdot \mid) : \mid
\cdot \mid \mbox{ is an equivalent norm on } X \}. $ 
\end{Def}

\noindent
We have that $\Delta (X) \ne \emptyset$ and it is easy to see that $\ddot{
\delta}_{\alpha}(X)= \sup \{ \gamma _{\alpha } : \gamma \in \ddot{\Delta}(X)\}$

\begin{Thm}\cite{OTW} \label{T:otw}
Let $X$ have a basis $(x_i)$. 
\begin{enumerate}
\item If $\gamma \in \Delta (X)$ then $\gamma _{\alpha }$ is a continuous 
decreasing 
function of $\alpha$. Also $\gamma _{\alpha + \beta} \geq
\gamma_{\alpha} \gamma_{\beta}$ for all $\alpha, \beta < \omega _1$.
\item For all $\alpha < \omega _1$ and $n \in \N$, $\ddot{\delta}_{\alpha \cdot
n}(X)=(\ddot{\delta}_{\alpha} (X))^n$.
\item $X$ does not contain $\ell_1$ iff $\ddot{\delta}_{\alpha}(X)=0$ for some
$\alpha < \omega _1$. 
\end{enumerate}
\end{Thm}

\begin{Def}
Let $X$ have a basis $(x_i)$. The {\em spectral index} $I_{\Delta}(X)$ is 
defined  by $I_{\Delta}(X) = \inf \{ \alpha < \omega _1 : \ddot{\delta}
_{\alpha}(X)<1 \}$ if such an $\alpha$ exists and $I_{\Delta}(X)=\omega_1$, 
otherwise. 
\end{Def}

\begin{Def}{\em Mixed Tsirelson Norms} \cite{AD}  
Let $F \subseteq \N$. Let $(\alpha _n)_{n \in F}$ be a set of countable 
ordinals and $(\theta _n)_{n \in F} \subset (0,1)$. The mixed Tsirelson space
$T(\theta_n , S_{\alpha _n})_{n \in F}$ is the completion of $c_{00}$ under
the implicit norm \[ \| x \| = \| x \| _{\infty} \vee \sup_{q \in \N} \sup \{
\theta_q \sum_1 ^n \| E_i x \| : (E_i)_1 ^n \mbox{ is an } \alpha _q 
\mbox{-admissible sequence of sets }\}. \] 
\end{Def}

\noindent
It is proved in \cite{AD} that such a norm exists.
They also proved that $T(\theta_n, S_{a_n})_{n \in F}$ is reflexive if $F$ is 
finite or $\lim_{F \ni n \rightarrow \infty} \theta_n =0$. $(e_n)$ is a 
1-unconditional basis for $T(\theta_n , S_{\alpha_n})$ so we can restrict the
$E_i$'s in the above definition to be {\em intervals}. It is worth noting 
that $T$, Tsirelson's space \cite{Ts} as described in 
\cite{FJ} satisfies $T=T(1/2, S_1)=T(1/2^n, S_n)_{\N}$.

\section{A property of the $\Delta$-spectrum} \label{S:Delta}

\noindent
The definition of $\delta_{\alpha}(x_i)$ is w.r.t. \ the coordinate
system  $(x_i)$. In \cite{OTW} the following notion is also introduced:

\begin{Def}
Let $(e_i)$ be a basis for $X$ and let $(x_i) \prec (e_i)$. For $\alpha < 
\omega_1$ we define 
\begin{eqnarray*}
\delta_{\alpha}((x_i),(e_i))= \sup \! \! \! & \{ & \! \! \! \delta \geq 0 : \| 
\sum _1 ^n  y_i \| \geq \delta \sum _1 ^n \| y_i \| \mbox{ whenever } 
(y_i)_1 ^n \prec (x_i) \\
& & \! \! \! \mbox{and } (y_i)_1 ^n \mbox{ is } \alpha \mbox{-admissible  
w.r.t. } (e_i) \}. 
\end{eqnarray*}
If $(y_i) \prec (e_i)$ we say that $(y_i)$ $\Delta 
_{(e_i)}$-stabilizes $\gamma =(\gamma_{\alpha})_{\alpha <\omega_1}$ if 
there exists $\varepsilon_n \searrow 0$ so that for
all $\alpha < \omega_1$ there exists $m \in \N$ so that if $n \geq m$ and 
$(z_i) \prec (y_i)_n ^{\infty}$ then $\mid \delta_{\alpha}((z_i),(e_i))- 
\gamma_{\alpha} \mid < \varepsilon_n$. Let $\Delta(X,(e_i))$ be the set of all
$\gamma$'s so that $(y_i)$ $\Delta_{(e_i)}$-stabilizes $\gamma$ for some
$(y_i) \prec (e_i)$. 
\end{Def}

\noindent
One can show, by the same arguments used to establish the analogous result
for $\Delta(X)$ \cite{OTW},  that for all $(x_i) \prec (e_i)$ there exists 
$(y_i) \prec (x_i)$
and $\gamma = ( \gamma_{\alpha})_{\alpha < \omega_1}$ so that $(y_i)$ $\Delta
_{(e_i)}$-stabilizes $\gamma$. In particular, $\Delta(X,(e_i))$ is 
non-empty. 

\bigskip
\noindent
In this section we prove that the $\Delta$-stabilization and the $\Delta_{(e_i)
}$-stabilization are actually the same notions. More precisely we prove

\begin{Thm} \label{T:stabil}
Let $X$ have a basis $(e_i)$ and let $(x_i) \prec (e_i)$ so that $(x_i)$ 
$\Delta_{(e_i)}$-stabilizes $\bar{\gamma} \in \Delta(X,(e_i))$ and $(x_i)$
$\Delta$-stabilizes $\gamma \in \Delta(X)$. Then $\bar{\gamma}=\gamma$.
Hence $\Delta(X)=\Delta(X,(e_i))$.
\end{Thm}

\noindent
First we need a combinatorial result. $[ \N ]$ denotes the set of infinite 
subsequences
of $\N$. If $N=(n_i) \in [ \N ]$ then $S_{\alpha}(N) = \{ (n_i)_{i \in F} : F
\in S_{\alpha} \}$ and $[N]$ is the set of infinite subsequences of $N$.

\begin{Prop} \label{P:salpha}
Let $N \in [ \N ]$. Then there exists $L=(\ell_i) \in [N]$ so that for 
all $\alpha < \omega_1$,
\[ (\ell_i)_{i \in F} \in S_{\alpha}
  \im (\ell_{i+1})_{i \in F} \in S_{\alpha}(N). \]
\end{Prop}

\noindent
{\bf Proof} Let $N= (n_i)$. We shall choose $M=(m_i) \in [ N ]$ and then
prove by induction on $\alpha$ that $L=(\ell_i)$ satisfies the proposition
where $\ell_i=n_{m_i}$.
Let $m_1 = n_1$. If $m_\subk$ has been defined set $m_{k+1}=n_{m_\subk}$. 

\bigskip
\noindent
The case $\alpha =0$ is trivial.

\bigskip
\noindent
Assume the result holds for $\alpha$ and that $(n_{m_i})_{i \in F} \in S_{
\alpha +1}$. Thus there exists $k \in \N$ and $n_{m_\subk} \leq E_1 < E_2 < \cdots
E_{n_{m_\subk}}$ (some possibly empty) so that $E_j \in S_{\alpha}$ for all $j$ and
$(n_{m_i})_{i \in F} = \cup_1 ^{n_{m_\subk}} E_j$. For each $j$ let $E_j=(n
_{m_i})_{i \in F_j}$. Then $n_{n_{m_{k}}}=n_{m_{k+1}} \leq (n_{m_{i+1}})_{i \in
F} = \cup_1 ^{n_{m_\subk}}(n_{m_{i+1}})_{i \in F_j}$ and for all $j$, $(n_{m_{i+1}
})_{i \in F_j} \in S_{\alpha}(N)$. Therefore $(n_{m_{i+1}})_{i \in F} \in S_{
\alpha +1}(N)$.

\bigskip
\noindent
If $\alpha$ is a limit ordinal and $\alpha _n \nearrow \alpha$ are the ordinals
used to define $S_{\alpha}$ and the result holds for all $\beta < \alpha$ (so
in particular for each $\alpha_n$), let $(n_{m_i})_{i \in F} \in S_{\alpha}$.
Thus for some $k \in \N$, $k \leq \min(n_{m_i})_{i \in F} \equiv n_{m_{i_0}} 
\leq (n_{m_i})_{i \in F} \in S_{\alpha _\subk}$. Hence $n_\subk \leq n_{n_{m_{i_0}}}=
n_{m_{i_0 +1}} \leq (n_{m_{i+1}})_{i \in F} \in S_{\alpha _\subk}(N)$ therefore
$(n_{m_{i+1}})_{i \in F} \in S_{\alpha}(N)$. \hfill $\Box$

\bigskip
\noindent
As a corollary we obtain a result of independent interest.

\begin{Cor}
Let $N \in [ \N ]$. Then there exists $L=(\ell_i) \in [ N ]$ so that for all
$\alpha < \omega_1$,
\[ (\ell_i)_{i \in F} \in S_{\alpha} \im (\ell_i)_{i \in F \backslash \min(F) }
  \in S_{\alpha}(N). \] 
\end{Cor}

\noindent
{\bf Proof} Let $L$ be as in proposition~\ref{P:salpha}. Let $F=(f_1 < f_2<
\cdots < f_r)$ with $(\ell_i)_{i \in F} \in S_{\alpha}$. Thus $(\ell_{f_1+1},
\ell_{f_2+1}, \ldots , \ell_{f_r+1}) \in S_{\alpha}(N)$. Since $f_1+1 \leq f_2$,
$f_2+1 \leq f_3, \ldots $ and $S_{\alpha}(N)$ is both spreading and 
hereditary we get that $(\ell_i)_{i \in F \backslash \min(F)} \in 
S_{\alpha}(N) $. \hfill $\Box$

\bigskip
\noindent
{\bf Proof of theorem~\ref{T:stabil}}
Let $(x_i)$ $\Delta_{(e_i)}$- and $\Delta$-stabilize $\bar{\gamma}$ and
$\gamma$ respectively and let $\alpha < \omega_1$. Since $S_{\alpha}$ is
spreading, $\bar{\gamma} \leq \gamma$. Let $\varepsilon >0$ and choose $(y_i)
\prec (x_i)$ so that for all $(z_i) \prec (y_i)$,
\[ \mid \delta_{\alpha}(z_i)- \gamma_{\alpha} \mid < \varepsilon. \]
For $i \in \N$ set $n_i = \min ( \ran(y_i))$ w.r.t. $(e_i)$ and choose $L=
(n_{m_i})$ by proposition~\ref{P:salpha}. For $w \in \span(y_{m_i})$ if 
$w=\sum_{i=j}^\ell a_i y_{m_i}$ where $a_j \ne 0$ we set $\bar{w}= \sum_{i=j+1}^\ell
a_i y_{m_i}$. 

\bigskip
\noindent
{\bf Claim:} If $(w_i)_1 ^\ell \prec (y_{m_i})$ is $\alpha$
-admissible w.r.t. $(e_j)$ then $(\bar{w}_i)_1^\ell$ is $\alpha$
-admissible w.r.t. $(y_j)$. 

\bigskip
\noindent
Indeed let  $m_{k_i}= \min ( \ran(w_i))$ w.r.t. $(y_j)$. Then $n_{m_{k_i}} =
\min ( \ran(w_i))$ w.r.t. $(e_j)$, and $(n_{m_{k_i}})_1^\ell \in S_{\alpha} \im
(n_{m_{k_i +1}})_1^\ell \in S_{\alpha}((n_j)) \im (m_{k_i+1})_1^\ell \in 
S_{\alpha}$. Since $m_{k_i +1} \leq \min (\ran( \bar{w}_i))$ w.r.t. $(y_j)$,
and $S_{\alpha}$ is spreading the claim follows.

\bigskip
\noindent
We may assume that $\| y_{m_i} \| =1$ for all $i$ and that no subsequence of 
$(y_{m_i})$ is equivalent to the unit vector basis of $c_0$ (indeed, if this
were false then clearly $\bar{\gamma}_0 = \gamma _0 =1$ and $\bar{\gamma}
_{\alpha} = \gamma _{\alpha} =0 $ for all $\alpha \geq 1$). Thus by taking
long averages of $(y_{m_i})$ we may choose $(z_i) \prec (y_{m_i})$ with the 
property that for all $z \in \span(z_i)$ 
\[ \| z-\bar{z} \| < \varepsilon \| \bar{z}\|. \] 
By the definition of $\bar{\delta}_{\alpha} \equiv 
\delta_{\alpha}((z_i),(e_i))$ there
exists $(w_i)_1^\ell \prec (z_i)$ which is $\alpha$-admissible w.r.t. $(e_j)$ and
satisfies 
\[ \| \sum _1^\ell w_i \| < (\bar{\delta}_{\alpha} + \varepsilon) \sum
_1^\ell \| w_i \|. \] 
By the above claim $(\bar{w}_i)_1^\ell$ is $\alpha$-admissible
w.r.t. $(y_j)$. Furthermore 
\begin{eqnarray*}
\| \sum_1^\ell \bar{w}_i \| & \leq & \| \sum _1^\ell w_i
\| + \sum_1^\ell \| w_i - \bar{w}_i \| < (\bar{\delta}_{\alpha}+ \varepsilon )
\sum_1^\ell \| w_i \| + \sum _1^\ell \varepsilon \| \bar{w}_i \| \\ 
& < & [(\bar{\delta}_{\alpha}+ \varepsilon)(1+ \varepsilon)+ \varepsilon ] 
\sum _1^\ell \| \bar{w}_i \|. 
\end{eqnarray*}
It follows that 
$\gamma_{\alpha} -\varepsilon < \delta_{\alpha}(y_i) < (\bar{\gamma}_{\alpha}
+ \varepsilon)(1+ \varepsilon ) + \varepsilon. $ Since $\varepsilon$ is
arbitrary we obtain $\gamma _{\alpha} \leq \bar{\gamma}_{\alpha}$ and so 
$\gamma_{\alpha} = \bar{\gamma}_{\alpha}$. 

\bigskip
\noindent
To  prove that $\Delta(X)=
\Delta(X,(e_i))$, let's first show the inclusion $\subseteq$. Let $(x_i)$ 
$\Delta$-stabilize $\gamma \in \Delta(X)$. We can find $(y_i) \prec
(x_i)$ that $\Delta_{(e_i)}$-stabilizes  $\bar{\gamma} \in \Delta_{(e_i)}$.
But then $(y_i)$ $\Delta$-stabilizes $\gamma$, therefore $\gamma= 
\bar{\gamma}$, thus $\gamma \in \Delta_{(e_i)}$. The inclusion $\supseteq$
is proved similarly.  \hfill $\Box$ 

\section{The space $T(\theta_n, S_n)_{\N}$} \label{S:mixedtsirelson} 

\noindent
If $\theta_n \not \rightarrow 0$ or if $\theta_n =1$ for some $n$ then 
$T(\theta_n,
S_n)_{\N}$ is isomorphic to $\ell_1$. Thus we shall confine ourselves to the
case where $\sup \theta_n <1$ and $\theta_n \rightarrow 0$. Furthermore we 
assume that $\theta_n \searrow 0$ and $\theta_{m+n}\geq \theta_n \theta_m$ for
all $n,m \in \N$. Indeed it is easy to see that $T(\theta_n, S_n)_{\N}$ is 
naturally isometric to $T(\bar{\theta}_n, S_n)_{\N}$ where 
\[ \bar{\theta}_n \equiv \sup
\{ \prod_{i=1}^\ell \theta_{k_i}: \sum_{i=1}^\ell k_i \geq n \}. \]

\begin{Def} A sequence $(\theta_n)$ of scalars is called regular if
$(\theta_n) \subset (0,1)$, $\theta_n \searrow 0$ and $\theta_{n+m} \geq \theta
_n \theta_m$ for all $n,m \in \N$. If the sequence $(\theta_n)$ is regular
we define the space $T(\theta_n, S_n)_{\N}$ to be regular. 
\end{Def}

\noindent
{\em Throughout this section, the spaces $T(\theta_n, S_n)_{\N}$ will always
be assumed to be regular}.

\bigskip
\noindent
It is easy to see (eg   \cite{OTW}) that if a sequence $(b_n) \subset (0,1]$ 
satisfies
$b_{n+m} \geq b_n b_m$ for all $n,m \in \N$ then $\lim_n b_n^{1/n}$ exists and
equals $\sup_n b_n^{1/n}$. Therefore, if the sequence $(\theta_n)$ is regular 
then the limit $\theta \equiv \lim_{n \rightarrow \infty} \theta_n ^{1/n} = 
\sup _n \theta_n^{1/n}$ exists. Note also that if $(X, \norm \cdot \norm)$ is a
Banach space with a basis, then $\delta_{n+m}(X) \geq \delta_n(X) \delta_m(X)$
for all $n,m \in \N$, thus $\lim_n \delta_n(X)^{1/n}= \sup_n \delta_n(X)^{1/n}$
exists. Furthermore, if $X$ does not contain $\ell_1$ isomorphically, then
$1> \delta_n(X) \searrow 0$.

\bigskip
\noindent
For $n \in \N$, define $\phi_n \equiv \frac{\theta _n}{\theta^n}$. We 
easily see
\begin{itemize}
\item If $\theta=1$ then $\phi_n = \theta_n \searrow 0$.
\item $\phi_{n+m} \geq \phi_n \phi_m$ for all $n,m \in \N$.
\item $\phi_n ^{1/n} \rightarrow 1$.
\item $\phi_n \leq 1, \forall n \in \N$.
\end{itemize}
>From now on, for a regular sequence $(\theta_n)$ we will be referring to the 
limit $\theta = \lim \theta_n^{1/n}$ and the representation $\theta_n = 
\theta^n \phi_n$ as above.

\bigskip
\noindent
The main theorem in this section is the following

\begin{Thm} \label{T:main}
Let $X=T(\theta_n, S_n)_{\N}$ be regular and let $\theta = \lim_n \theta_n 
^{1/n}$. Then
\begin{itemize}
\item[(1)] For all $Y \prec X$, $\ddot{\delta}_1(Y)=\theta$. Moreover for all 
$\varepsilon >0$ there exists an equivalent norm $\mid \cdot \mid$ on $X$ so
that $\delta_1((X,\mid \cdot \mid),(e_i)) > \theta - \varepsilon$. 
\item[(2)] For all $Y \prec X$ and for all $n \in \N$, $\ddot{\delta}_n(Y)=
 \theta^n$
    and $\ddot{\delta}_{\omega}(Y)=0$.
\item[(3)] $ \mbox{ For all } Y \prec X,  I_{\Delta}(Y) = \left\{ \begin{array}{ll}
    \omega & \mbox{ if } \theta =1 \\
    1 & \mbox{ if } \theta < 1  
    \end{array} \right. $
\item[(4)] If $\frac{\theta_n}{\theta^n} \rightarrow 0$ then $X$ is arbitrarily
 distortable.
\end{itemize}
\end{Thm}

\noindent
To prove the above theorem we need the following two results

\begin{Prop} \label{P:delta1geqtheta}
Let $X=T(\theta_n, S_n)_{\N}$ be regular. Then for every $\varepsilon >0$ 
there is an equivalent 1-unconditional norm $\mid \cdot \mid$ on $X$ such 
that $\delta_1 ((X,\mid \cdot \mid), (e_i)) \geq \theta - \varepsilon$.
\end{Prop}

\begin{Thm} \label{T:computedeltaj}
Let $X=T(\theta_n,S_n)_{\N}$ be regular. Then for all $Y \prec X$ and $j \in \N
$ we have \[ \delta_j(Y) \leq \theta^j \sup_{p \geq j} \phi _p \vee \frac{
\theta_j}{\theta_1}. \]
\end{Thm}

\bigskip
\noindent
{\bf Proof of theorem~\ref{T:main}}\\ 
{\bf (1)} To prove that if $Y \prec X$ then $\ddot{\delta}_1(Y) \leq \theta$ 
we note that 
if $\norm \cdot \norm$ is an equivalent norm on $T(\theta_i, S_i)_{\N}$
then there exists $C \geq 1$ such that $C^{-1} \delta_n(Y)
\leq \delta_n(Y, \norm \cdot \norm) \leq
C \delta_n(Y)$ for all $n \in \N$. Let $\delta_n \equiv \delta_n(Y,
\norm \cdot \norm)$. Then since for all $n$ and  $m$, $\delta_{n+m} \geq 
\delta_n
\delta_m$ we have $\lim_n \delta_n^{1/n}= \sup_n \delta_n^{1/n}$ exists.
Hence $\delta_1 \leq \lim \delta_n^{1/n}= \lim \delta_n(Y)^{1/n}$, the latter
limit existing for the same reason. Now 
\[ \lim \delta_n(Y)^{1/n} \leq \lim_{n \rightarrow \infty} (\theta^n 
\sup_{p \geq n}\phi_p \vee \frac{\theta_n}{\theta_1})^{1/n} = \theta \]
by theorem~\ref{T:computedeltaj}. Thus $\ddot{\delta}_1(Y) \leq \theta$
as was to be proved.
The ``moreover'' part is proposition~\ref{P:delta1geqtheta} and this completes 
the proof of $\ddot{\delta}_1(Y)= \theta$.

\noindent
{\bf (2)} Since $\ddot{\delta}_1(Y)= \theta$ we obtain $\ddot{\delta}_n(Y)=
\theta^n$ from theorem~\ref{T:otw}. By theorem~\ref{T:computedeltaj} we
have that for all $\gamma \in \Delta(Y)$ and for all $j \in \N$, 
\[ \gamma_j \leq \theta^j \sup_{p \geq j} \phi_p \vee \frac{\theta_j}{
\theta_1}. \]
Therefore, again by theorem~\ref{T:otw}, for all $\gamma \in \Delta(Y)$, 
$\gamma_{\omega} = \lim _{n \in \N} \gamma_n =0$. Hence, for every equivalent
norm $\mid \cdot \mid$ on $Y$, for every $\gamma \in \Delta (Y, \mid \cdot 
\mid)$, $\gamma_{\omega}=0$. Since $\ddot{\delta}_{\omega}(Y)= \sup \{ \gamma
_{\omega}: \gamma \in \ddot{\Delta}(Y) \}$ we have $\ddot{\delta}_{\omega}(Y)=
0$. 

\noindent
{\bf (3)} Follows immediately from (2)

\noindent
{\bf (4)} Let $\lambda>1$. Choose $n \in \N$ so that $\sup_{p \geq n}\phi_p < 
\frac{\theta_1}{2 \lambda}$.
By (1) we can define an equivalent norm $\norm \cdot \norm $ on $X$
such that 
\[ \delta_n((X,\norm \cdot \norm),(e_i))
\geq \delta_1 ((X,\norm \cdot \norm),(e_i))^n \geq \frac{\theta^n}{2}. \]
Let $Y \prec X$. 
By  equation (\ref{E:equiv}) of section~\ref{S:prel},
there exists $C>0$ such that 
\[  C \norm y \norm \leq \| y \| \leq  D(Y, \norm \cdot \norm) C \norm y
 \norm, \mbox{ for all } y \in  Y. \]
Therefore by Observation~\ref{O:equiv}, 
\[ D(Y, \norm \cdot \norm)  \geq 
\frac{\delta_n(Y, \norm \cdot \norm)}{\delta_n(Y, \| \cdot \|)}. \]
Since  $\delta_n(Y, \norm \cdot \norm) \geq \delta_n((X,
\norm \cdot \norm),(e_i)) \geq \frac{\theta^n}{2}$, and $\delta_n(Y, \| \cdot 
\|) \leq \theta^n \sup_{p \geq n}\phi_p \vee \frac{\theta_n}{\theta_1} \leq
\frac{1}{\theta_1} \theta^n \sup _{p \geq n} \phi_p$ (by
theorem~\ref{T:computedeltaj}),
we obtain $D(Y, \norm \cdot \norm)  \geq \frac{\theta_1}{2 \sup_{p \geq n} 
\phi_p}>\lambda$. \hfill $\Box$

\bigskip
\noindent
The proof of proposition~\ref{P:delta1geqtheta} comes from an argument in 
\cite{OTW}. We recall this argument here.

\bigskip
\noindent
{\bf Sketch of the proof of proposition~\ref{P:delta1geqtheta}}     
Fix $n \in \N$ such that $\theta_n^{1/n} > \theta - \varepsilon$ and set 
$a \equiv \theta_n ^{1/n}$ For $j \in \N$ and $x \in X$ define
\begin{eqnarray*}
& & \hskip.6in \mid x \mid_j = \sup \{ a^j \sum_1^\ell \| E_i x \| : (E_i x)_1
  ^\ell \mbox{ is } j \mbox{-admissible w.r.t. } (e_i) \} \mbox{ and }  \\
& & \hskip.6in \mid x \mid =\frac{1}{n} \sum_{j=0}^{n-1} \mid x \mid_j  
  (\mbox{where } \mid \cdot \mid_0 = \| \cdot \|).
\end{eqnarray*} 
We claim that $\delta_1((X, \mid 
\cdot \mid ), (e_i)) \geq a$. To see this let $e_\subk \leq x_1 < x_2 < \cdots
< x_\subk$ in $X$ and $x= \sum_{i=1}^k x_i$. 
For $ j=1, \ldots , n-1$ we have
$ \mid x \mid _j \geq a \sum_{i=1}^k \mid x_i \mid_{j-1} $ (by the definitions
of $\mid \cdot \mid_j$ and $\mid \cdot \mid_{j-1}$) 
 and also $ \mid x \mid_0 \geq a \sum_{i=1}^k \mid x_i \mid _{n-1}$ (since
$a^n= \theta_n$). Therefore we get $ \mid x \mid \geq a \sum_{i=1}^k 
\mid x_i \mid $. \hfill $\Box$

\bigskip
\noindent
To prove theorem~\ref{T:computedeltaj} we need some norm estimates in
$T(\theta_n,S_n)_{\N}$ for certain iterated rapidly increasing averages. Before
defining what we mean by this we fix some terminology.

\bigskip
\noindent
Let $E$ be an interval in $\N$ and $x \in c_{00}$. We say 
that {\em $E$ does not split $x$} \/ if either $E \cap \ran(x)= \emptyset$ or 
$\ran(x) \subseteq E$. 
Let $(x_i)$ be a block basis of $(e_i)$ in $c_{00}$, $x \in \span(x_i)_i$, 
$N \in \N$, and $E_1<E_2< \cdots <E_N$ be intervals in $\N$ so that $\cup_{i
=1}^N E_i \subseteq \ran(x)$. We say that we {\em minimally shrink the 
intervals $(E_\subl)_{\ell =1}^N$ to obtain intervals $(F_\subl)_{\ell =1}^n$ 
which don't split the $x_i$'s}, if for $\ell=1, \ldots ,N$ we let $G_\subl =
E_\subl \backslash \cup \{ \ran(x_i): E_\subl \mbox{ splits } x_i \}$ and let
$F_1<F_2< \cdots F_n$ be the enumeration of the non-empty $G_\subl$'s.

\bigskip
\noindent
By a {\em  tree} \/ we shall mean a non-empty partially ordered  set 
$({\cal T}, \ll)$ for which the set $\{ y \in {\cal T}:y \ll x \}$ is linearly
 ordered and finite for each $x \in {\cal T}$. If ${\cal T}' \subseteq 
{\cal T}$ then we say that $({\cal T}', \ll)$ is a {\em subtree} \/ of $({\cal
T},\ll)$. The tree ${\cal T}$ is called
{\em finite} \/ if the set ${\cal T}$ is finite. The {\em initial} \/  nodes of 
${\cal T}$ are the minimal elements of ${\cal T}$ and the {\em terminal} 
\/ nodes are
the maximal elements. A {\em branch} \/ in ${\cal T}$ is a maximal linearly 
ordered set
in  ${\cal T}$. The {\em immediate successors} \/ of $x \in {\cal T}$
are all  the
nodes $y \in {\cal T}$ such that $x \ll y$ but there is no $z \in {\cal T}$ 
with $x \ll z \ll y$.
If $X$ is a linear space, then a {\em tree in $X$} \/ is a tree 
whose nodes are vectors in $X$. 
If $X$ is a Banach space with a basis $(e_i)$ and $(x_i) \prec (e_i)$
then an {\em admissible  averaging tree of $(x_i)$}, is a finite 
tree ${\cal T}$ in $X$ with the following properties:
\begin{itemize}
\item ${\cal T} = (x^j_i)_{j=0,i=1}^{M,N^j}$ where
  $M \in \N$ and $1=N^M \leq \cdots \leq N^1 \leq N^0$.
\item $x^j_1 < \cdots <x^j_{N^j}$ w.r.t. $(e_s)$ ($j=0,1, \ldots ,
  M-1$) $\&$ $(x^0_i)_{i=1}^{N^0}$ is a subsequence of $(x_s)$.
\end{itemize}
Also for $j=1, \ldots , M$ and $i=1, \ldots , N^j$ we have the
following:
\begin{itemize}
\item There exists a non-empty interval $I^j_i \subseteq \{1, \ldots ,
  N^{j-1} \}$ such that $\{ x^{j-1}_s : s \in I^j_i \}$ are the
  immediate successors of $x^j_i$.
\item $x^j_i = \frac{1}{ \mid I^j_i \mid }\sum_{s \in I^j_i}
  x^{j-1}_s$.
\item $(\min(\ran(x^{j-1}_s)))_{s \in I^j_i} \in S_1$ where
  $\ran(x^{j-1}_s)$ is taken w.r.t. $(x_s)$.
\end{itemize}
Note that the last two properties together require that $x^j_i$ be a
1-admissible average of all of its immediate successors w.r.t. $(x_s)$.
Let ${\cal T}= (x^j_i)_{j=0,i=1}^{M,N^j}$ be an
admissible averaging tree as in the above definition, and let $b= \{
y_M \ll \cdots \ll y_0 \}$ be a branch in ${\cal T}$. For $i=0, 1,
\ldots , M $ we say that the {\em level} \/ of $y_i$ is $i$. Note that
this is well defined, since the definition of admissible averaging
trees forces every branch to have the same number of elements. Indeed
for each $i$ and $j$, the level of $x^j_i$ in ${\cal T}$ is $j$.
Let ${\cal T}$ be a tree, $x \in {\cal T}$ \/
of level $\ell$ and $k \in
\N$. By ${\cal T}(x, k)$ (resp. ${\cal T}^*(x,k)$) we shall denote the subtree
of ${\cal T}'=\{ x \} \cup \{ y \in {\cal T}: y \gg x \}$ (resp. 
${\cal T}'= \{ y \in {\cal T}: y \gg x \}$) that contains all the nodes of 
${\cal T}'$ that have level $\ell, \ell -1, \ldots , \mbox{ or } \ell -k +1$ in
${\cal T}$.
Let ${\cal T}$ be an admissible averaging tree in a Banach 
space $X$ with a basis $(e_i)$, $x \in {\cal T}$  with immediate
successors $x_1< \cdots <x_n$ (a finite block basis of $(e_i)$), $k
\in \N$,  and let $F 
\subseteq \N$ be an interval which does not split any of $x_1, \ldots
, x_n$.  Then by 
${\cal T}_F(x,k)$ we shall denote the subtree of ${\cal T}(x,k)$ given by 
${\cal T}_F(x,k)=\{ x \} \cup \{ y \in {\cal T}^*(x,k): \ran(y) \mbox{
  (w.r.t. $(e_i)$)} \subseteq F \}$.

\begin{Def} \label{D:meaverage}
Let $(x_i)$ be a block sequence of $(e_i)$ in $c_{00}$, $M, N \in \N$, and let
$(\varepsilon^j_i)_{j,i \in \N} \subset (0, 1)$.  We say that $x$ is 
an $(M, (\varepsilon^j_i), N)$ average of $(x_i)$ w.r.t. $(e_i)$ if 
there exists an admissible averaging tree ${\cal T}= (x^j_i)_{j=0,i=1}^{M,
N^j}$ of $(x_i)$  whose initial node is $x (=x^M_1)$ and 
\begin{quote}
for $j=1, \ldots, M$ and $1 \leq i \leq N^j$ if $N^j_i= \max(\ran(x^j_i))$ w.r.t
$(e_s)$ ($N^j_0=N$), then $x^j_i$ is an average of its immediate successors
of length $k^j_i > \frac{2 N^j_{i-1}}{\varepsilon^j_i}$.
\end{quote}
${\cal T}$ then will be called an $(M,(\varepsilon^j_i),N)$ admissible 
averaging tree of $(x_i)$ w.r.t. $(e_i)$. For $i=1, \ldots , N^0$ set 
$N^0_i = \max(\ran(x^0_i))$ w.r.t. $(e_s)$, and $N^0_0=N$. Then $(N^j_i)
_{j=0, i=0}^{M,N^j}$ are called the maximum coordinates of ${\cal T}$ w.r.t.
 $(e_i)$.
\end{Def}

\begin{Rmk} \label{R:meaverage}
Let $X$ be a Banach space with basis $(e_i)$ and let $(x_i)$ be a block 
sequence 
of $(e_i)$ with $\| x_i \| \leq 1$ for all $i \in \N$. Let  $(\varepsilon^j_i)
_{j,i\in \N} \subset (0,1)$. Let  
$M, N \in \N$ and let $x$ be an $(M, (\varepsilon^j_i), N)$ average of 
$(x_i)$ w.r.t. $(e_i)$ given by ${\cal  T}=(x^j_i)_{j=0,i=1}^{M,N^j}$. Then  
we can write $x= \sum_{i \in F} a_i x_i$ for some finite set $F
\subset  \N$ such that 
\begin{itemize}  
\item[(1)] $\sum_{i \in F}a_i=1$ $ \& $ $a_i > 0$ for all $i \in F$. 
\item[(2)] $x$ is $M$-admissible w.r.t. $(x_i)$ (i.e. $F \in S_M$).
\item[(3)] Let $(N^j_i)_{j=0,i=0}^{M,N^j}$ be the 
  maximum coordinates of ${\cal T}$ w.r.t. $(e_s)$. For $j=1, \ldots ,
  M$  and $1   \leq i \leq  N^j$,  let 
  $E^j_i(1)<E^j_i(2)< \cdots < E^j_i(N^j_{i-1})$ be a finite sequence of 
  intervals in $\N$ with   $\cup_{\ell=1}^{N^j_{i-1}} E^j_i(\ell) \subseteq 
  \ran(x^j_i)$ and assume that we minimally 
  shrink the $E^j_i(\ell)$'s to obtain intervals  $(F^j_i(\ell))_{\ell=1}^
  {N^j_{i-1}}$ (some of which may be empty) which don't split the $x^{j-1}_i$'s.
  Then 
  \[ \sum_{j=1}^M \sum_{i=1}^{N^j} \sum_{\ell =1}^{N^j_{i-1}} \| (E^j_i( \ell)
  \backslash F^j_i(\ell))x^j_i \| < \sum_{j,i} \varepsilon^j_i. \]
\end{itemize}
\end{Rmk}

\noindent
Indeed {\em (1)} and {\em (2)}  are obvious. To see {\em (3)} 
note  that  for every $j=1, \ldots , M$, $1 \leq i \leq N^j$ 
and $\ell =1, \ldots , N^j_{i-1}$,  the set $E^j_i(\ell)$ splits at most two
$x^{j-1}_s$'s each of them having norm at most 1. Thus  
$\| (E^j_i(\ell) \backslash F^j_i(\ell))x^j_i \| \leq 2/k^j_i$ and so 
$\sum_{\ell=1}^{N^j_{i-1}} \| E^j_i( \ell) \backslash F^j_i(\ell) \|< 2 N^j_{i-1
}/k^j_i < \varepsilon^j_i$, which proves  {\em (3)}.

\bigskip
\noindent
The concept of $(M, (\varepsilon^j_i), N)$ vectors is implicit in \cite{AD} 
(see also \cite{OTW}).

\begin{Prop} \label{P:meaverage}
Let $(x_i)$ be a block sequence in $c_{00}$, $M,N \in \N$ and $(\varepsilon^j
_i)_{j,i \in \N} \subset (0,1)$. Then there exists $x$ which is an $(M,
(\varepsilon^j_i), N)$ average of $(x_i)$ w.r.t. $(e_i)$. 
\end{Prop}

\noindent
{\bf Proof}
Note that by replacing each $(\varepsilon^j_i)_i$ by a smaller sequence if
necessary we may assume that $(\varepsilon^j_i)_i$ is decreasing.
For $M=1$ we choose $x^1_1$ to be an average of $k^1_1>2N / \varepsilon^1
_1$ many $x_s$'s chosen from $\{ x_s :s \geq k^1_1 \}$. Next, consider the case 
$M=2$. At first we continue the argument that we gave for $M=1$ to construct
$\bar{x}^1_1 < \bar{x}^1_2< \cdots$ as follows: For $\bar{k}^1_1> 2N/ 
\varepsilon^1_1$ let $\bar{x}^1_1$ be an average of $\bar{k}^1_1$ many $x_s$'s
chosen from $\{ x_s : s \geq \bar{k}^1_1 \}$. 
If $\bar{x}^1_i$ has been constructed for some $i \in \N$, and $\bar{k}^1_{i+1}
> 2\bar{N}^1_i / 
\varepsilon^1_{i+1}$, then $\bar{x}^1_{i+1}$ is taken to be an average
of  $\bar{k}^1_{i+1}$
 many $x_s$'s chosen from $\{ x_s : s \geq \bar{k}^1_{i+1} \}$ where $\bar{N}^1
_i = \max(\ran(\bar{x}^1_i))$ w.r.t. $(e_s)$. Note
that $\bar{x}^1_i < \bar{x}^1_{i+1}$ since $\varepsilon^1_{i+1} < 1$. Also note
 that for
every $i \in \N$, $\bar{x}^1_i$ is a 1-admissible w.r.t. $(x_s)$.
Then for $k^2_1 > 2N / \varepsilon^2_1$ take $x^2_1$ to be an average of
$k^2_1$ many $\bar{x}^1_s$'s chosen from $\{ \bar{x}^1_s : \bar{x}^1_s \geq 
x_{k^2_1} \}$. Then
the $(2, (\varepsilon^j_i),N)$ admissible averaging tree ${\cal T}$ of $(x_i)$
that corresponds to $x^2_1$ is determined as follows: $x^2_1 \in {\cal T}$. If 
$x^2_1= \frac{1}{ \mid F \mid} \sum_{i \in F} \bar{x}^1_i$ for some finite set 
$F 
\subset \N$ then $\bar{x}^1_i \in {\cal T}$ for $i \in F$. For each $i \in F$ 
if $\bar{x}^1
_i =\frac{1}{\mid F_i \mid} \sum _{s \in F_i}x_s$ for some finite set $F_i 
\subset \N$ then $x_s \in {\cal T}$ for $s \in F_i$. Enumerate the $x_s$'s in 
${\cal T}$ as $x^0_1<x^0_2< \cdots <x^0_{N^0}$ and the $\bar{x}^1_s$'s in 
${\cal T}$ 
as $x^1_1< x^1_2< \cdots <x^1_{N^1}$. Since $x^2_1$ is a 1-admissible average 
of $(x^1_i)$ w.r.t. $(x_i)$ and for each $i=1, \ldots , N^1$, $x^1_i$ is 
1-admissible w.r.t.  $(x_s)$,  we have that $x^2_1$ is 2-admissible w.r.t. 
$(x_i)$. We let the $k^1_i$'s and $N^1_i$'s be defined by 
definition~\ref{D:meaverage}. Each $k^1_i$ will be $\bar{k}^1_{i'}$
for Some $i' \geq i$
and $N^1_0=N$ while $N^1_i=N^1_{i'}$. Since $(\varepsilon^1_i)$ is decreasing
the condition $k^1_i > 2N^1_{i-1} / \varepsilon^1_i$ remains valid. The case 
$M>2$ is proved by iterating this procedure.  \hfill $\Box$ 

\begin{Rmk} \label{R:refine}
Definition~\ref{D:meaverage} requires only that $k^j_i > 2N^j_{i-1} / 
\varepsilon^j
_i$ The proof shows that we could also construct $(x^j_i)$ so that $k^j_i >
\frac{6 (N^j_{i-1})^2}{ \theta_1 \varepsilon^j_i}$. We will use this remark
in lemma~\ref{L:cbsequences}.
\end{Rmk}

\bigskip
\noindent
Next we  prove some norm estimates for $(M, (\varepsilon^j_i) ,N)$ averages in 
$T(\theta_n,S_n)_{\N}$. $\| \cdot \|$ will always denote the norm of $T(\theta
_n, S_n)_{\N}$. We need for $p \in \N \cup \{ 0 \}$ and $N \in \N$ to define 
the equivalent norms $\| \cdot \|_p$ and  $\| \cdot \|_{S_N,p}$ and the 
continuous seminorms 
$\| \cdot \|_{N, p}$ as follows ($\| \cdot \|_0 = \| \cdot \|$  and $\theta_0
=1$):
\begin{eqnarray*}
& & \hskip.5in 
\| x \|_p  = \theta_p \sup \{ \sum \| E_i x \| :(E_i) \mbox{ is a } p 
                \mbox{-admissible sequence of intervals } \} \\
& & \hskip.5in 
\| x \|_{N,p} = \sup \{ \sum_1^N \| E_i x \|_p : N \leq E_1 <E_2< 
        \cdots < E_N \mbox{ are intervals }\} \mbox{ and }\\
& &  \hskip.5in 
\| x \|_{S_N,p} =\sup \{ \sum \| E_i x \|_p :(E_i) \mbox{ is an } N 
         \mbox{-admissible sequence of intervals } \}.
\end{eqnarray*}      
Of course for $x \in c_{00}$ each ``sup'' above is a ``max''
and there exists $p \in \N$ so that $\| x \| = \| x \|_p$ if $\| x \| \not =
\| x \|_{\infty}$.

\begin{Rmk} \label{R:ratiotheta}
Let $\theta_0=1$. For all $x \in c_{00}$ and for all $p \in \N$ we have
\[ \| x \| _p \leq \frac{\theta_p}{\theta_{p-1}} \| x \|_{S_1,p-1}.\]
Moreover if $p=1$ we have equality.
\end{Rmk}

\noindent
Indeed there exists $(E_i)_{i \in I}$ a $p$-admissible family of intervals such
that 
\[ \| x \|_p = \theta_p \sum_{i \in I} \| E_i x \|. \]
We can write $I=\cup_1^\ell I_j$ where $ (E_i)_{i \in I_j}$ is $p-1$-admissible 
and if $F_j$ is the smallest interval including $\cup_{i \in I_j}E_i $ then
$(F_j)_1^\ell$ is 1-admissible. Thus 
\[ \| x \|_p = \frac{\theta_p}{\theta_{p-1}}
\sum_{j=1}^\ell \theta_{p-1} \sum_{i \in I_j} \| E_i x \| \leq \frac{\theta
_p}{\theta_{p-1}} \sum_{j=1}^\ell \| F_j x \|_{p-1} \leq
\frac{\theta_p}{\theta_{p-1}} \| x \|_{S_1,p-1}.  \hskip.5in \Box  \]
{\bf Notation} If $A \subset [0, \infty)$ is a finite non-empty set, we set
$A^*=A \backslash \{ \max(A) \}$. 

\begin{Obs} \label{O:star}
Let $N \in \N$ and $D, \varepsilon >0$. Note
that if $k \geq \frac{ND}{\varepsilon}$ and $A_\subl \subset [0,D]$ for $\ell
=1, 
\ldots ,N$ are finite sets with $\mid A_1 \mid + \cdots + \mid A_N \mid \leq k$
then $\frac{1}{k} \sum_{\ell=1}^N \sum \{ a: a \in A_\subl \} \leq \max( \cup
_{\ell=1}^N A^*_\subl) + \varepsilon$. 
\end{Obs}

\noindent
We will apply this for $D=\frac{1}{\theta_1}$ in the proof of {\it (2)}\/ of 
lemma~\ref{L:longaverage} below. 

\begin{Lem} \label{L:longaverage} 
Let non-zero vectors $ k \leq x_1<x_2<\ldots <x_\subk$ with $\| x_i\| \leq 1$ 
for all $i$, $x=\frac{1}{k}(x_1 + \cdots + x_\subk)$ and $\varepsilon
\in (0, 1)$. 
Let $F \subseteq \ran(x)$ be an interval in $\N$ which does not split the 
$x_i$'s. Set $\theta_0 =1$, $N_i = \max ( \ran(x_i))$ w.r.t. $(e_j)$, $N_0=1$
 and let $N \in \N$. If $k > \frac{6 N}{\theta_1 \varepsilon}$ then 
\begin{itemize}
\item[(1)] For every $p \in \N$, 
$ \| Fx \|_{N,p} \leq \frac{\theta_p}{\theta_{p-1}} \max \{ 
  \| x_i \|_{N_{i-1},p-1} : \ran(x_i)  \subseteq F \} + \varepsilon$.   
\item[(2)] There exists $n \in \N$, intervals  $F_1<F_2< \ldots <F_n$ which 
don't split any $x_i$, $\cup_{\ell =1}^n F_\subl \subseteq \ran(x)$, and 
$(p_\subl)_{\ell =1}^n \subset \N$ so that 
\[ \| x \|_{N,0} \leq \max \left(\bigcup_{\ell=1}^n \{ \frac{\theta_{p_\subl}}{
\theta_{p_\subl-1}} \| x_i \|_{N_{i-1},p_\subl-1}: \ran(x_i) \subset F_\subl
 \}^* \right) + \varepsilon. \]
\end{itemize}
\end{Lem}

\noindent
{\bf Proof}
{\bf (1)} For $p \in \N$ there exist intervals $N \leq E_1< \ldots <E_N$ such 
that $\cup_{\ell =1}^N E_\subl \subseteq F$ and
\[ \| Fx \|_{N,p}  =\sum_{\ell =1}^N \| E_\subl x \|_p \leq \frac{\theta_p}{
\theta_{p-1}}  \sum_{\ell=1}^N \| E_\subl x \|_{S_1,p-1} \mbox{ (by 
Remark~\ref{R:ratiotheta}) }. \] 
We minimally shrink the intervals $(E_i)_1^N$ to get $n \leq N$ and intervals
$N \leq F_1<F_2< \cdots <F_n$ which don't split the $x_i$'s.
Since each $E_\subl$ splits at most two $x_i$'s, $\| \cdot \|_{S_1,p-1} \leq
\frac{1}{\theta_1} \| \cdot \|$ and $\frac{\theta_p}{\theta_{p-1}} \leq 1$, 
\[ \frac{\theta_p}{\theta_{p-1}} \sum_{\ell=1}^N \| E_\subl x \|_{S_1,p-1} \leq 
\frac{\theta_p}{\theta_{p-1}} \sum_{\ell=1}^n \| F_\subl x \|_{S_1,p-1} +
\frac{2N}{k \theta_1}. \]
Fix an $\ell \in \{ 1, \ldots , n \}$. 
There exists a 1-admissible family of 
intervals $(F_{\ell,m})_m$  with $F_{\ell,m} \subseteq F_\subl$ for all 
 $m$ and $ \| F_\subl x \|_{S_1,p-1}  = \sum_m \| F_{\ell,m}x \|_{p-1}$. 
Let $s$ be minimal
with $\ran(x_s) \cap F_{\ell,1}\not = \emptyset$ (we may assume that such an 
$s$ exists) and $t$ be maximal with $\ran(x_t) \cap F_\subl \not = \emptyset$. 
Then
\begin{eqnarray*}
\sum_m \| F_{\ell,m}x \|_{p-1} & \leq & \frac{1}{k} ( \sum_m \| F_{\ell ,m}x_s \| + 
  \| x_{s+1} \|_{N_s, p-1}+ \cdots + \| x_t \|_{N_s,p-1} ) \\
   & \leq & \frac{1}{k} ( \frac{1}{\theta_1} + \sum_{i=s+1}^t \| x_i \|
   _{N_{i-1},p-1}).
\end{eqnarray*}
Set $[r,R] \equiv \{ i : \ran(x_i) \subseteq F \}$. Hence
\[  \frac{\theta_p}{\theta_{p-1}}\sum_{\ell=1}^n \| F_\subl x \|_{S_1,p-1} 
 \leq  \frac{\theta_p}{\theta_{p-1}} \frac{1}{k} (\frac{n}{\theta_1}+
 \| x_{r+1} \|_{N_r,p-1}+ \| x_{r+2} \|_{N_{r+1},p-1}+\cdots + \| x_R \|_
{N_{R-1},p-1}).\] 
Therefore we have proved that
\[ \| Fx \|_{N,p} \leq \frac{1}{k} \frac{\theta_p}{\theta_{p-1}} \left( \| x_{r
+1} \|_{N_r,p-1} + \| x_{r+2} \|_{N_{r+1},p-1} + \cdots + \| x_R \|_{N_{R-1},
p-1} \right)  + \frac{3N}{k \theta_1}. \]
Thus
\begin{equation} \label{E:1/2}
\| Fx \|_{N,p} \leq \frac{1}{k} \frac{\theta_p}{\theta_{p
-1}} \sum \{ \| x_i \| _{N_{i-1},p-1}: \ran(x_i) \subseteq F \} + \frac{3N}{k 
\theta_1}.
\end{equation}
This yields (1).

\bigskip
\noindent
{\bf (2)} Choose intervals $N \leq E_1<E_2< \ldots <E_N$ such that 
$\| x \|_{N,0}= \sum_{l=1}^N \| E_\subl x \|$. As before, we minimally shrink  
the intervals $(E_i)$ to obtain $n \leq N$ and non-empty intervals $F_1<
F_2< \cdots <F_n$ which don't split the $x_i$'s and satisfy
\[ \sum_{\ell =1}^N \| E_\subl x \| \leq \sum_{\ell =1}^n \| F_\subl x \| + 
\frac{2N}{k}. \] 
Fix $\ell \in \{ 1, \ldots , n\}$. If $\| F_\subl x \| \not = \|
F_\subl x  \|_\infty$ there exists 
$p_\subl \in \N$ such that $\| F_\subl x \| = \| F_\subl x \|_{p_\subl}$.
By equation (\ref{E:1/2}) for $N=1$ we get
\begin{equation} \label{E:3/4}
\| F_\subl x \|_{p_\subl}  \leq  \frac{1}{k} \frac{\theta_{p_\subl}}{\theta
_{p_\subl-1}} \left( \sum \{ \| x_i \|_{N_{i-1},p_\subl-1}: 
\ran(x_i) \subseteq F_\subl \} \right) +\frac{3}{k \theta_1}. 
\end{equation}
If $\| F _\subl x \| = \| F_\subl x \|_{\infty}$ then $\| F_\subl x \|
\leq \frac{1}{k}$ and so (\ref{E:3/4}) still is valid. Thus
\begin{eqnarray*}
\| x \|_{N,0} & \leq &\frac{1}{k} \sum_{\ell=1}^n 
\frac{\theta_{p_\subl}}{\theta_{p_\subl-1}} \sum \{ 
\| x_i \|_{N_{i-1},p_\subl-1}: \ran(x_i) \subseteq F_\subl \} + \frac{5 N}{k 
\theta_1} \\
 & < & \max \left(\bigcup_{\ell =1}^n \{ \frac{\theta_{p_\subl}}{
\theta_{p_\subl-1}} \| x_i \|_{N_{i-1},p_\subl-1}: \ran(x_i) \subset F_\subl
 \}^* \right) + \varepsilon 
\end{eqnarray*}
by observation~\ref{O:star} since $\| \cdot \|_{N_{i-1}, p_\subl -1} \leq 
\frac{1}{\theta_1} \| \cdot \|$, and $k > \frac{6 N}{\varepsilon \theta_1} =
\frac{N}{(\varepsilon /6) \theta_1}$. \hfill $\Box$

\bigskip
\noindent
Combining lemma~\ref{L:longaverage} with proposition~\ref{P:meaverage} and 
remark~\ref{R:refine} we obtain

\begin{Lem} \label{L:cbsequences}
Let $(x_i)$ be a normalized block sequence in $X=T(\theta_i,S_i)_{\N}$, $ M,N
 \in \N$ and  $(\varepsilon^j_i)_{j, i \in \N} \subset (0,1)$.
There exists $x$, an $(M, (\varepsilon^j_i),N)$ average of $(x_i)$ w.r.t. 
$(e_i)$, so that if ${\cal T}=(x^j_i)_{j=0,i=1}^{M,N^j}$ is the $(M, 
(\varepsilon^j_i),N)$ admissible averaging tree of $(x_i)$ with $x=x^M_1$,
and $(N^j_i)_{j=0,i=0}^{M,N^j}$ are the maximum coordinates of ${\cal T}$ 
w.r.t. $(e_i)$ then for $j=1, \ldots , M$ and $i=1, \ldots , N^j$ we have the 
following properties:
\begin{itemize}
\item[(1)] For every $p \in \N$ and every $F \subseteq \ran(x^j_i)$ which does 
  not split any $x^{j-1}_s$ we have
\[ \| Fx^j_i \|_{N^j_{i-1},p} 
 \leq\frac{\theta_p}{\theta_{p-1}}    \max \{ \| x^{j-1}_s\|
    _{N^{j-1}_{s-1},p-1}: \ran(x^{j-1}_s) \subseteq F \} +\varepsilon^j_i /
  N^j_{i-1}. \]
\item[(2)] There exists $n \in \N$ and  intervals $F_1<F_2<\ldots <F_n$
  which don't split any $x^{j-1}_s$, ($\cup_{\ell=1}^n F_\subl \subseteq \ran(x
  ^j_i$))  and $(p_\subl)_{\ell =1}^n \subseteq \N$   such that 
\[ \| x^j_i \|_{N^j_{i-1} , 0} \leq \max \left( \bigcup_{\ell=1}^n \{  
 \frac{\theta_{p_\subl}}{\theta_{p_\subl-1}}  \| x^{j-1}_s \|_{N^{j-1} _{s-1},
 p_\subl-1} : ran( x^{j-1}_s ) \subseteq F_\subl \}^* \right) + 
\varepsilon^j_i.\]
\end{itemize}
\end{Lem}  

\begin{Lem} \label{L:cbsequences345}
Let $(x_i)$ be a normalized block sequence in $X=T(\theta_i , S_i)_{\N}$,
$\varepsilon >0$, $(\varepsilon^j_i)_{j,i \in \N} \subset (0,1)$ with
$\sum_{j,i} \varepsilon^j_i < \varepsilon$ and let $x$ be an $(M, 
(\varepsilon^j_i),N)$ average of $(x_i)$ w.r.t. $(e_i)$. Let ${\cal T}= (x^j_i)
_{j=0,i=1}^{M,N^j}$ be the $(M,(\varepsilon^j_i),N)$ admissible averaging tree 
of $(x_i)$ with $x=x^M_1$, let $(N^j_i)_{j=0,i=0}^{M,N^j}$ be the 
maximum
coordinates of ${\cal T}$ w.r.t. $(e_i)$ and assume that for $j=1, \ldots , M$ 
and $i=1,\ldots , N^j$ the properties (1) and (2) of lemma~\ref{L:cbsequences} 
are satisfied. Then we have
\begin{itemize}
\item[(3)] If $0 \leq p' < p$, $p-p' \leq j \leq M$, $1 \leq i \leq N^j$
  and $F \subseteq \ran(x^j_i)$ is an interval which does not split any $x^{j-1
 }_s$ then 
\[ \| Fx^j_i \|_{N^j_{i-1},p} \leq \frac{\theta_p}{\theta_{p'}}\max  \{ 
 \| x^{j-(p-p')}_s \|_{N^{j-(p-p')}_{s-1},p'}:
 \ran( x^{j-(p-p')}_s) \subseteq F \}  + 
 \sum_{x^k_s \in {\cal T}_F(x^j_i,p-p') }  \frac{\varepsilon^k_s}{N^k_{s-1}}. \]
\item[(4)] If $1 \leq p \leq j \leq M$, $1 \leq i \leq N^j$ and $F\subseteq \N$ 
 is an interval which does not split any $x^{j-1}_s$ then
 \[ \| Fx^j_i \|_{N^j_{i-1},p} \leq   \theta_p \max \{ \| x^{j-p}_s 
    \|_{N^{j-p}_{s-1},0} : \ran(x^{j-p}_s) \subseteq F \} +
   \sum \{ \frac{\varepsilon ^k_s}{N^k_{s-1}} :x^k_s \in {\cal T}_F (x^j_i,p) 
   \}. \]
\item[(5)] There exists $m \in \N$ and intervals $F_1<F_2< 
  \ldots <F_m$ ($\cup_\subl F_\subl \subseteq \ran(x^M_1)$) which don't
  split the $x^0_s$'s and $(p_\subl)_{\ell=1}^m \subset \N$ with $p_\subl \geq 
 M$  for all  $\ell$, such that 
\[\| x^M_1 \| \leq \max \left( \bigcup_{\ell=1}^m \{ \frac{\theta
 _{p_\subl}}{\theta _{p_\subl-M}} \| x^0_s \|_{N^0_{s-1}, p_\subl-M}:
 \ran(x^0_s) \subseteq F_\subl \}^* \right) + \varepsilon. \] 
\end{itemize}
\end{Lem}

\noindent
{\bf Proof} \\
{\bf (3)} By (1) of lemma~\ref{L:cbsequences} we have
\begin{eqnarray*} & & \| Fx^j_i \|_{N^j_{i-1},p}  \leq  \frac{\theta_p}{
\theta_{p-1}}  \max \{ \| x^{j-1}_s \|_{N^{j-1}
_{s-1},p-1}: \ran(x^{j-1}_s)\subseteq F \} + \frac{\varepsilon^j_i}{N^j_{i-1}}\\
 & \leq & \frac{\theta_p}{\theta_{p-1}} \frac{\theta_{p-1}}{\theta_{p-2}}
  \max \{ \| x^{j-2}_s \|_{N^{j-2}_{s-1},  p-2}: \ran(x^{j-2}_s) \subseteq F \}
  +  \sum \{\frac{ \varepsilon^k_s}{N^k_{s-1}} : x^k_s \in {\cal T}_F (x^j_i,2)
    \} \\ 
 & \leq & \ldots \\
 & \leq & \frac{\theta_p}{\theta_{p-1}} \frac{\theta_{p-1}}{\theta_{p-2}} 
  \ldots \frac{\theta_{p'+1}}{\theta_{p'}} \max  \{ 
  \| x^{j-(p-p')}_s \|_{N^{j-(p-p')}_{s-1},p'}: \\
 & & \hskip 1.7in  \ran(x^{j-(p- p')}_s) \subseteq F \} +  \sum \{ \frac{
   \varepsilon^k_s}{N^k_{s-1}} :x^k_s \in {\cal T}_F (x^j_i, p-p') \}. 
\end{eqnarray*}
{\bf (4)} Follows immediately from (3), letting $p' =0$.\\
{\bf (5)} We prove by induction on $J$ that 
\begin{quote}
for  $J=1, \ldots , M$ and $1 \leq i \leq N^J$
there exists $m \in \N$, intervals $F_1<F_2< \cdots <F_m$ ($\cup _\subl F_\subl
 \subseteq \ran(x^J_i)$) that don't split the $x^0_s$'s, and $(p_\subl)_{\ell=
 1}^m \subset \N$ with $p_\subl \geq J$ for all $\ell$, such that 
\[ \| x^J_i \|_{N^J_{i-1},0} \leq \max \left( \bigcup_{\ell=1}^m \{ \frac{
\theta  _{p_\subl}}{\theta _{p_\subl-J}} \| x^0_s \|_{N^0_{s-1}, p_\subl-J}:
 \ran(x^0_s) \subseteq F_\subl \}^* \right) + 
\sum \{ \varepsilon^k_s : x^k_s \in {\cal T}( x^J_i,J) \} \]
\end{quote}  
((5) then follows by taking $(J,i)=(M,1)$ and noting that $\| x^M_1 \| \leq
\| x^M_1 \|_{N,0}=\|x \|_{N^M_0,0}$). Indeed, for $J=1$ this follows from
the statement of (2) for $j=1$. Assume that the statement is proved for all 
positive integers $\leq J$ 
where $J\leq M-1$. By (2) there exist intervals $F_1'< \cdots
<F_n'$ ($\cup_\subl F_\subl ' \subseteq \ran(x^{J+1}_i)$) which don't split the
 $x^J_s$'s,  and $(p_\subl ')_{\ell=1}^n$ such that  
\[ \| x^{J+1}_i \|_{N^{J+1}_{i-1},0} \leq \max \left( \bigcup_{\ell=1}
 ^n \{ \frac{\theta_{p_\subl '}}{\theta_{p_\subl '-1}} \| x^J_s \|_{N^J_{s-1},
 p_\subl '-1}: \ran(x^J_s) \subseteq F_\subl ' \}^* \right) + 
 \varepsilon ^{J+1}_i. \]
If $p_\subl '-1 =0$ for some $\ell$ and $\ran(x^J_s) \subseteq F_\subl '$ then
by the induction hypothesis there exists $M(s) \in \N$, intervals $F_1(s)<F_2(s)
< \cdots < F_{M(s)}(s)$ $(\cup_ \mu F_\mu (s) \subseteq \ran(x^J_s))$ that don't
split the $x^0_t$'s and $(p_\mu (s))_{\mu =1}^{M(s)} \subset \N$ with $p_\mu (s)
\geq J$ for all $\mu$ such that
\[ \| x^J_s \|_{N^J_{s-1},0} \leq \max \left( \bigcup_{\mu =1}^{M(s)} \{
\frac{\theta_{p_\mu (s)}}{\theta_{p_\mu (s) -J}} \| x^0_t \|_{N^0_{t-1},p_\mu 
(s)-J} : \ran(x^0_t) \subseteq F_\mu (s) \}^* \right)+ \sum \{ \varepsilon^k_t 
 : x^k_t \in {\cal T}(x^J_s,J) \}. \]
If $0< p_\subl '-1 \leq J$ for some $\ell$, and $\ran(x^J_s) \subseteq F_\subl 
'$  then by (4), 
\[ \| x^J_s \|_{N^J_{s-1},p_\subl '-1} \leq \theta_{p_\subl '-1} \max \{ \| x^{
J-p_\subl '+1}_t  \|_{N^{J-p_\subl '+1}_{t-1},0}: \ran(x^{J-p_\subl '+1}_t) 
\subseteq \ran(x^J_s) \} + \sum_{x^k_t \in {\cal T}(x^J_s, p_\subl ' -1) } 
\varepsilon^k_t. \]
For the remaining $\ell$'s we have by (3) for $j=J$, $p=p_\subl ' -1$ and 
$p'=p_\subl ' -1-J$,
\[ \| x^J_s \|_{N^J_{s-1},p_\subl ' -1} \leq \frac{\theta_{p_\subl ' -1}}{
\theta_ {p_\subl ' -1 -J }} \max
 \{ \| x^0_t \|_{N^0_{t-1},p_\subl '-1-J}: \ran(x^0_t) \subseteq \ran(x^J_s) 
 \} + \sum \{ \varepsilon^k_t : x^k_t \in {\cal T}(x^J_s,J) \}. \]
Combining  these estimates we get 
\begin{eqnarray*}
& & \| x^{J+1}_i \|_{N^{J+1}_{i-1},0} \leq  \\
& & \max \left( \bigcup_{\{ \ell : p_\subl '=1 \} } \bigcup_{ \{ s: \ran(x^J_s)
  \subseteq F_\subl ' \} } \bigcup_{ \mu =1}^{M(s)} \{ \frac{ \theta_1
  \theta _{p_\mu (s)}}{\theta_{p_\mu (s)-J}} \| x^0_t \|_{N^0_{t-1},
  p_\mu(s)-J} + \sum
  \{ \varepsilon^k_w : x^k_w \in {\cal T}(x^J_s,J) \} : \ran(x^0_t) \subseteq
  F_\mu(s) \}^* \right. \\
& & \cup \bigcup_{\{ \ell : 0< p_\subl '-1 \leq J \}} \{ \theta_{p_\subl '}  
  \| x^{J-p_\subl '+1}_t \| 
   _{N^{J-p_\subl '+1}_{t-1},0}+ \sum \{ \varepsilon^k_s : x^k_s \in {\cal T}^*
  (x^{J+1}_i, p_\subl ') \} : 
   \ran(x^{J-p_\subl '+1}_t) \subseteq F_\subl ' \}^*  \\
& & \left. \cup \bigcup_{ \{ \ell : p_\subl '>J+1 \} } \{ \frac{
   \theta_{p_\subl '}}{\theta_{p_\subl '-(J+1)}} \| x^0_t \|_{N^0
    _{t-1},p_\subl '- (J+1)} + \sum \{ \varepsilon^k_s : x^k_s \in {\cal T}^*
  (x^{J+1}_i, J+1) \} : \ran(x^0_t) \subseteq F_\subl ' \}^* \right) 
     + \varepsilon^{J+1}_i.
\end{eqnarray*}
The induction hypothesis gives that for $0< p_\subl '-1 < J$ and $1 \leq t \leq
N^{J-p_\subl '+1}$ with $\ran(x^{J-p_\subl '+1}_t) \subseteq F_\subl '$, there 
exists $K(\ell ,t) \in \N$
 and sets $G_1(\ell ,t)<G_2(\ell , t)< \ldots < G_{K(\ell , t)}(\ell ,t)$ which
 don't split the $x^0_s$'s such that $\cup_\subk G_\subk (\ell , t) \subseteq  
\ran(x^{J-p_\subl '+1}_t)$, and there exist $(q_\subk (\ell , t)) _{k =1}^{K(
\ell, t)} \subset \N$ with $q_\subk (\ell , t) \geq J-p_\subl '+1$ such that 
\begin{eqnarray*}
\| x^{J-p_\subl '+1}_t \|_{N^{J-p_\subl '+1}_{t-1},0} & \leq & \max \left( 
\bigcup_{k=1}^{K(\ell , t)} \{ 
 \frac{\theta_{q_\subk(\ell ,t)}}{\theta_{q_\subk(\ell , t)-(J-p_\subl '+1)}} 
\| x^0_s \|_{N^0_{s-1},  q_\subk(\ell ,t)-(J-p_\subl '+1)}:  \ran(x^0_s) 
\subseteq G_\subk (\ell ,t)\}^* \right) \\
& + &   \sum \{  \varepsilon^k_s : x^k_s \in {\cal S}(\ell ,t)  \}
\end{eqnarray*}
where ${\cal S}(\ell ,t)={\cal T}(x^{J-p_\subl '+1}_t, J-p_\subl +1)$. Thus, 
these estimates give
\begin{eqnarray*}
& & \| x^{J+1}_i \|_{N^{J+1}_{i-1},0} \leq  \\
& & \max \left( \bigcup_{\{ \ell : p_\subl '=1 \} } \bigcup_{ \{ s: \ran(x^J_s)
  \subseteq F_\subl ' \} } \bigcup_{\mu =1}^{M(s)} \{ \frac{\theta_1
  \theta_{p_\mu (s)}}{\theta_{(1+ p_\mu (s))-(J+1)}} \| x^0_t
  \|_{N^0_{t-1}, p_\mu(s)-J}  \right. \\
& & \hskip 1.3in + \sum \{ \varepsilon^k_w : x^k_w \in {\cal T}(x^J_s,J) \} 
  : \ran(x^0_t)   \subseteq   F_\mu(s) \}^* \\
& & \cup \bigcup_{ \{ \ell :0< p_\subl '-1<J \}} \bigcup_{k=1}^{K(\ell ,t)} \{ 
 \frac{ \theta_{p_\subl '} \theta_{q_\subk (\ell ,t)}}{\theta_{(p_\subl '+
  q_\subk (\ell ,t))-(J+1)}} \| x^0_s \|_{N^0_{s-1},  q_\subk(\ell ,t)-(J-
  p_\subl '+1)} \\
& & \hskip 1.3in  + \sum \{\varepsilon^k_s : x^k_s \in {\cal T}^*(x^{J+1}_i, 
  p_\subl ') \cup {\cal S}(\ell ,t) \} : \ran(x^0_s) \subseteq G_\subk (\ell ,
  t) \}^* \\
& & \cup \bigcup_{\{ \ell : p_\subl '=J+1 \}} \{ \theta_{J+1}  
  \| x^0_t \|_{N^0_{t-1},0}+ \sum \{ \varepsilon^k_s : x^k_s \in {\cal T}^*
  (x^{J+1}_i, J+1) \} : \ran(x^0_t) \subseteq F_\subl ' \}^*  \\
& & \left. \cup \bigcup_{ \{ \ell : p_\subl '>J+1 \} } \{ \frac{
   \theta_{p_\subl '}}{\theta_{p_\subl '-(J+1)}} \| x^0_t \|_{N^0
    _{t-1},p_\subl '- (J+1)} + \sum \{ \varepsilon^k_s : x^k_s \in {\cal T}^*
  (x^{J+1}_i, J +1) \} : \ran(x^0_t) \subseteq F_\subl ' \}^* \right) 
     + \varepsilon^{J+1}_i.
\end{eqnarray*}
Note that $\theta_1 \theta_{p_\mu (s)} \leq \theta_{1+p_\mu (s)}$, $1+p_\mu(s)
\geq J+1$, $\theta_{p_\subl '} \theta_{q_\subk(\ell ,t)} \leq \theta_{p_\subl '+
q_\subk (\ell ,t)}$, $p_\subl '+q_\subk (\ell ,t) \geq p_\subl '
+(J-p_\subl '+1)=J+1$, the sets $F_\mu (s)$'s $F_\subl' $'s and $ G_\subk(\ell 
,t)$'s don't 
split the $x^0_s$'s, and arranged in successive order, give the required 
sequence $F_1< \ldots < F_m$. Then $1+p_\mu(s)$'s, $p_\subl ' + q_k (\ell ,t)
$'s, and $p_\subl '$'s for $p_\subl ' \geq J+1$ arranged in the corresponding
order, give the required sequence $(p_\subl)_{\ell =1}^m$. This finishes the 
induction. \hfill $\Box$

\bigskip
\noindent
Combining lemmas~\ref{L:cbsequences} and~\ref{L:cbsequences345} we immediately
obtain

\begin{Cor} \label{C:normmeaverage}
Let $(x_i)$ be a normalized block sequence in $X=T(\theta_i,S_i)_{\N}$, $M,N 
\in \N$, $\varepsilon >0$ and  $(\varepsilon^j_i)_{j, i \in \N} \subset (0,1)$
 with $\sum_{j,i} \varepsilon^j_i < \varepsilon $. 
There exists $x$ an $(M, (\varepsilon^j_i), N)$
average of $(x^0_i)$ w.r.t. $(e_i)$, so that if ${\cal T}= (x^j_i)_{j=0,i=1}
^{M,N^j}$ is the admissible averaging tree of $(x_i)$ with $x=x^M_1$, 
and $(N^j_i)_{j=0,i=0}^{M,N^j}$ are the maximum coordinates of ${\cal T}$ 
w.r.t. $(e_i)$, then 
\begin{itemize}
\item[(1)]  For $j=1, \ldots , M$, $i=1, \ldots , N^j$, $1 \leq p \leq j $ and 
 an interval $F \subseteq \ran(x^j_i)$ which does not split any $x^{j-1}_s$,
 \[ \| Fx^j_i \|_{N^j_{i-1},p} \leq   \theta_p \max \{ \| x^{j-p}_s 
    \|_{N^{j-p}_{s-1},0} : \ran(x^{j-p}_s) \subseteq F \} +
   \sum \{ \frac{\varepsilon^k_s}{N^k_{s-1}} : x^k_s \in {\cal T}_F (x^j_i, p) 
  \}. \]
\item[(2)] There exists $m \in \N$ and intervals $F_1<F_2<   \ldots <F_m$ which
 don't  split the $x^0_s$'s and $(p_\subl)_{\ell=1}^m \subset \N$ with $p_\subl 
\geq M$ for all  $\ell$, such that 
\[\| x^M_1 \| \leq \max \left( \bigcup_{\ell=1}^m \{ \frac{\theta
 _{p_\subl}}{\theta _{p_\subl-M}} \| x^0_s \|_{N^0_{s-1}, p_\subl-M}:
 \ran(x^0_s) \subseteq F_\subl \}^* \right) + \varepsilon. \] 
\end{itemize}
\end{Cor}

\bigskip
\noindent
To prove theorem~\ref{T:computedeltaj} we need also the following 

\begin{Lem} \label{L:phi}
For all $J,N \in \N$, $\varepsilon >0$ and $Y \prec X=T(\theta_i , S_i)_\N$
there exists $y \in Y$ with $\|y\|=1$ and 
\[ \| y \|_{N,p} < \phi _p(1+ \varepsilon), \mbox{ for all } p=1, \ldots J. \]
\end{Lem}

\noindent
{\bf Proof} If this were false, then $ \exists J, N \in \N$ $\exists 
\varepsilon \in (0, 1/2)$ $\exists Y \prec X$ such that 
\begin{equation} \label{E:phi}
\| y \| \leq \max_{1 \leq p \leq J} \frac{1}{\phi_p (1+ \varepsilon)}  
   \| y \|_{N,p} \mbox{ for all } y \in Y. 
\end{equation}
Since $(1+ \varepsilon)^n 
\phi _{J(n+1)} \rightarrow \infty$ as $n \rightarrow \infty$ we may  choose $n 
\in \N$ such that
\[ 1> \frac{1}{(1+ \varepsilon)^n \theta^J \theta_1 \phi_{J(n+1)}}+
2\varepsilon. \]
Let $(x_s)$ be a normalized block sequence in $Y$ and apply 
corollary~\ref{C:normmeaverage} to $(x_s)$ for $(M, \varepsilon,N)=(J(n+1), 
\varepsilon \theta_{J(n+1)}\phi_1^J,N)$ for an appropriate sequence 
$(\varepsilon^j_i)$, to construct $x=\sum a_s x^0_s$,
a $(J(n+1),(\varepsilon^j_i),N)$ average of $(x_s)$ w.r.t. $(e_s)$. Let $x$ 
have a corresponding 
admissible averaging tree ${\cal T}=(x^j_i)_{j=0,i=1}^{J(n+1),N^j}$, 
and let the maximum coordinates of ${\cal T}$ be $(N^j_i)_{j=0,i=0}^{J(n+1),
N^j}$ w.r.t. $(e_i)$. Define  $\delta^j_i=\varepsilon^j_i /
\phi_1^J$ for $j,i \in \N$ and note that $\sum \delta^j_i < \varepsilon 
\theta_{J(n+1)}$. Note that if $1 \leq p \leq J$ then $\phi_p \geq \phi_1^p \geq
\phi_1^J$ and if $k,s \in \N$ then we have that $\frac{\varepsilon^k_s}{\phi_p}
\leq \delta^k_s$. There exists $1 \leq p^1 \leq J$ so that 
\[ \theta_{J(n+1)} = \theta_{J(n+1)} \sum \| a_s x^0_s \| \leq \| x \| \leq
\frac{1}{ \phi_{p^1}(1+\varepsilon)}   \| x  \|_{N,p^1} = \frac{1}{ 
\phi_{p^1}(1+\varepsilon)} \| x^{J(n+1)}_1 \|_{N^{J(n+1)}_0,p^1} \]
(since $N = N^{J(n+1)}_0$). Then by corollary~\ref{C:normmeaverage} (1), 
there exists $s^1 \in \N$ so that 
$\ran(x^{J(n+1)-p^1}_{s^1}) \subseteq \ran(x^{J(n+1)}_1)$ and also there 
exists a family of intervals $(E_i)_{i=1,\ldots , N^{J(n+1)-p^1}_{s^1-1}}
\subseteq \ran(x^{J(n+1)}_1)$ so that 
\begin{eqnarray*}
\theta_{J(n+1)}& \leq & \frac{1}{\phi_{p^1}(1+\varepsilon)}  \left( \theta_{p^1}
   \| x^{J(n+1)-p^1}_{s^1} \|_{N^{J(n+1)-p^1}_{s^1-1},0} 
 +\sum \{ \varepsilon^k_s: x^k_s \in {\cal T}(x^{J(n+1)}_1, p^1) \} \right)\\  
& \leq & \sum_{i=1}^{N^{J(n+1)-p^1}_{s^1-1}} \frac{\theta^{p^1}}{1+ 
 \varepsilon}   \| E_i x^{J(n+1)-p^1} _{s^1} \| 
 +\sum \{ \delta^k_s: x^k_s \in {\cal T}(x^{J(n+1)}_1, p^1) \}
\end{eqnarray*}
We minimally shrink  the $E_i$'s if necessary, to obtain $(F_i)$ which don't 
split the $x^{J(n+1)-p^1-1}_s$'s. Let ${\cal A}$ be the set of $x^{J(n+1)-p^1
-1}_s$'s that is split by the $E_i$'s. Thus we get
\[ \theta_{J(n+1)} \leq \sum_i \frac{\theta^{p^1}}{1+ \varepsilon} 
   \| F_i x^{J(n+1)-p^1}_{s^1} \| +2 \sum \{ \| x^k_s \| : x^k_s \in {\cal A}
 \} +\sum \{ \delta^k_s: x^k_s \in {\cal T}(x^{J(n+1)}_1, p^1) \}. \]
Similarly by~\ref{E:phi} for each $i$  there exist $1\leq p^2_i \leq J$ so that
(note that $N \leq N^{J(n+1)-p^1}_{s^1-1}$)
\begin{eqnarray*}
& & \theta_{J(n+1)} \leq\\
& & \sum_i \frac{\theta^{p^1}}{(1+ \varepsilon)} 
               \frac{1}{\phi_{p^2_i}(1 +\varepsilon)} 
   \| F_i x^{J(n+1)-p^1}_{s^1} \|_{N,p^2_i} + 2 \sum \{ \| x^k_s \| : x^k_s \in
  {\cal A} \} + \sum \{ \delta^k_s: x^k_s \in {\cal T}(x^{J(n+1)}_1, p^1) \} 
   \\
& \leq & \sum_i \frac{\theta^{p^1}}{(1+ \varepsilon)} 
               \frac{1}{\phi_{p^2_i}(1 +\varepsilon)} 
   \| F_i x^{J(n+1)-p^1}_{s^1} \|_{N^{J(n+1)-p^1}_{s^1-1},p^2_i} 
   +2 \sum \{ \| x^k_s \| : x^k_s \in  {\cal A} \}
 +\sum _{x^k_s \in {\cal T}(x^{J(n+1)}_1, p^1) } \delta^k_s.
\end{eqnarray*}
Then by corollary~\ref{C:normmeaverage} (1), for each $i$ there exists 
$s_i^2 \in \N$ so that $\ran(x^{J(n+1)-p^1-p_i^2}_{s_i^2}) \subseteq F_i$ and 
also family of intervals $(E_{i,j})_{j=1, \ldots ,N^{J(n+1)-p^1-p^2_i}
   _{s_i^2-1}} \subseteq \ran(x^{J(n+1)-p^1}_{s^1} )$ so that
\begin{eqnarray*}
& & \theta_{J(n+1)} \leq \\
& & \sum_i \frac{\theta^{p^1}}{(1+ \varepsilon)} 
    \frac{1}{\phi_{p^2_i}(1 + \varepsilon)} \left( \theta_{p^2_i}
  \| x^{J(n+1)-p^1-p^2_i}_{s^2_i} \|_{N^{J(n+1)-p^1-p^2_i}_{s^2_i-1},0} 
   + \sum \{ \frac{\varepsilon^k_s}{N^k_{s-1}} : x^k_s \in {\cal T}_{F_i}
   (x^{J(n+1)-p^1}_{s^1}, p^2_i)  \} \right)   \\
& & +2 \sum \{ \| x^k_s \| : x^k_s \in  {\cal A} \} +
 \sum \{ \delta^k_s: x^k_s \in {\cal T}(x^{J(n+1)}_1, p^1) \} \leq \\
& & \sum_i \sum_{j=1}^{N^{J(n+1)-p^1-p^2_i}_{s^2_i-1}} 
   \frac{\theta^{p^1+p^2_i}}{(1+\varepsilon)^2} 
   \| E_{i,j} x^{J(n+1)-p^1-p^2_i}_{s^2_i} \| +2 \sum\{ \| x^k_s \| : x^k_s \in 
  {\cal A} \}  + \sum \{ \delta^k_s : x^k_s \in {\cal S} \}
\end{eqnarray*}
where ${\cal S}={\cal T}(x^{J(n+1)}_1,p^1) \cup \cup_i \cup \{ {\cal T} 
(x^{J(n+1)-p^1-1}_t,p^2_i) : \ran(x^{J(n+1)-p^1-1}_t) \subseteq F_i \}$.
We increase ${\cal A}$ by including every node $x^{J(n+1)-p^1-p^2_i-1}_s$ which
is split by some $E_{i,j}$ and 
minimally shrink the $E_{i,j}$'s to get intervals $(F_{i,j})$ which don't split
the $x^{J(n+1)-p^1-p^2_i-1}_s$'s. Thus
\[ \theta_{J(n+1)} \leq \sum_i \sum_j \frac{\theta^{p^1+p^2_i}}{(1+ 
\varepsilon)^2} 
\| F_{i,j} x^{J(n+1)-p^1-p^2_i}_{s^2_i} \| +2 \sum \{ \| x^k_s \| : x^k_s \in  
{\cal A} \} + \sum \{ \delta^k_s : x^k_s \in {\cal S} \}. \]
For every $i,j$  there exists $1 \leq p^3_{i,j} \leq J$ so that we have (note
also that $N \leq N^{J(n+1)-p^1-p^2_i}_{s^2_i-1}$)
\begin{eqnarray*}
& & \theta_{J(n+1)} \leq \\
& & \sum_i \sum_j \frac{\theta^{p^1+p^2_i}}{(1+ \varepsilon)^2} 
\frac{1}{\phi_{p^3_{i,j}}(1+ \varepsilon)}
\| F_{i,j} x^{J(n+1)-p^1-p^2_i}_{s^2_i} \|_{N,p^3_{i,j}} +2 \sum \{ \| x^k_s \| 
: x^k_s \in {\cal A} \} + \sum \{ \delta^k_s : x^k_s \in {\cal S} \} \\
& \leq & \sum_i \sum_j \frac{\theta^{p^1+p^2_i}}{(1+ \varepsilon)^2} 
\frac{1}{\phi_{p^3_{i,j}}(1+ \varepsilon)}
\| F_{i,j} x^{J(n+1)-p^1-p^2_i}_{s^2_i} \|_{N^{J(n+1)-p^1-p^2_i}_{s^2_i-1} ,
p^3_{i,j}} + 2 \sum \{ \| x^k_s \| : x^k_s \in {\cal A} \} \\
& & + \sum \{ \delta^k_s : x^k_s \in {\cal S} \}.
\end{eqnarray*}
By corollary~\ref{C:normmeaverage} (1), for each $i,j$ there exists $s^3_{i,j} 
\in \N$ so that
\begin{eqnarray*}
& & \theta_{J(n+1)} \leq \sum_i \sum_j \frac{\theta^{p^1+p^2_i}}{(1+ 
\varepsilon)^2} 
\frac{1}{\phi_{p^3_{i,j}}(1+ \varepsilon)} \left( \theta_{p^3_{i,j}}
\| x^{J(n+1)-p^1-p^2_i-p^3_{i,j}}_{s^3_{i,j}} \|_{N^{J(n+1)-p^1-p^2_i-p^3_{i,
j}}_{s^3_{i,j}-1} ,0} \right.\\
& & \left. + \sum \{ \frac{\varepsilon^k_s}{N^k_{s-1}} : x^k_s \in {\cal T}
_{F_{i,j}}(x^{J(n+1)- p^1-p^2_i}_{s^2_i },p^3_{i,j}) \} \right) 
+2 \sum \{ \| x^k_s \| : x^k_s \in {\cal A} \} + 
\sum \{ \delta^k_s : x^k_s \in {\cal S} \} \\
& & \leq \sum_i \sum_j \frac{\theta^{p^1+p^2_i+p^3_{i,j}}}{(1+ \varepsilon)^3} 
\| x^{J(n+1)-p^1-p^2_i-p^3_{i,j}}_{s^3_{i,j}} \|_{N^{J(n+1)-p^1-p^2_i-p^3_{i,
j}}_{s^3_{i,j}-1} ,0} + 2 \sum \{ \| x^k_s \| : x^k_s \in {\cal A} \} + 
\sum \{ \delta^k_s : x^k_s \in {\cal S'} \}
\end{eqnarray*}
for some ${\cal S}' \subseteq {\cal T}$. We continue passing to lower levels
of the tree until we obtain $Jn \leq p^1 + p^2_i + \cdots + p^r_{i, \ldots ,k}
\leq J(n+1)-1$. On each branch of the tree 
we stop when this is satisfied. Thus we get an estimate 
of the following form (${\cal A}$ increases to contain the $x^k_s$'s that are 
split)
\begin{eqnarray*}
\theta_{J(n+1)} & \leq & \sum \sum_{i, \ldots , k} 
  \frac{\theta^{p^1+p^2_i+ \cdots +p^r_\subk}}{(1+ \varepsilon)^r} 
 \| x^{J(n+1)-p^1-p^2_i \cdots -p^r_{i, \ldots , k}}_{s^r_{i,\ldots ,k}} \|
 _{N^{J(n+1)-p^1-p^2_i \cdots -p^r_{i, \ldots , k}}_{s^r_{i,\ldots ,k}-1},0} \\ 
& & + 2 \sum \{ \| x^k_s \| : x^k_s \in {\cal A} \} + 
\sum \{ \delta^k_s : x^k_s \in {\cal W} \}  
\end{eqnarray*}
for some ${\cal W} \subseteq {\cal T}$,
where the first ``$\sum $'' is taken over all branches on which we have
$Jn \leq p^1 + p^2_i + \cdots + p^r_{i, \ldots ,k} \leq J(n+1)-1$.
By remark~\ref{R:meaverage} (3) we have
that 2 $\sum \{ \| x^k_s \| : x^k_s \in {\cal A} \}   < \varepsilon 
\theta_{J(n+1)}$. Also $\sum \{ \delta ^k_s : x^k_s \in {\cal W} \} < 
\varepsilon \theta_{J(n+1)}$. Thus 
\[ \theta_{J(n+1)} \leq  \frac{\theta^{Jn}}{(1+ \varepsilon)^n} \sum \sum_{i, 
\ldots , k}  \| x^{J(n+1)-p^1-p^2_i \cdots - p^r_{i, \ldots , k}}_{s^r_{i,
\ldots ,k}} \|_{N^{J(n+1)-p^1-p^2_i \cdots -p^r_{i, \ldots , k}}_{s^r_{i,\ldots
 ,k}-1},0}  + 2 \varepsilon \theta_{J(n+1)}. \] 
Since $\| \cdot \|_{n,0} \leq \frac{1}{\theta_1} \| \cdot \|$, the vectors 
$x^{J(n+1)-p^1-p^2_i \cdots - p^r_{i, \ldots , k}}_{s^r_{i,\ldots ,k}}$'s have
disjoint support and their level in the tree is at least 1,by the triangle 
inequality we obtain
\[ \theta_{J(n+1)} \leq \frac{\theta^{Jn}}{\theta_1 (1+ 
\varepsilon)^n} + 2 \varepsilon \theta_{J(n+1)} \eq 
 1 \leq \frac{1}{(1+ \varepsilon)^n \theta^J \theta_1 \phi_{J(n+1)}} + 
  2 \varepsilon \] which is a contradiction. \hfill $\Box$ 

\bigskip
\noindent
{\bf Proof of theorem~\ref{T:computedeltaj}} Let $\varepsilon >0$ be
arbitrary. By lemma~\ref{L:phi} we can find a normalized block
sequence $(x_i)$ in $Y$ and an increasing  sequence $(\bar{j}_i)$ of integers,
$\bar{j}_1=1$,  so that if $N_0=1$ and $N_i = \max( \ran(x_i))$ w.r.t. 
$(e_s)$ then for every $i \in \N$ we have
\[ \begin{array}{ll} 
 \forall p=1, \ldots , \bar{j}_i, & \| x_i \|_{N_{i-1},p}< 
        \phi_p(1 + \varepsilon) \mbox{ and } \\
        \forall p \geq \bar{j}_{i+1}, & \| x_i \|_{N_{i-1},p}< \varepsilon. 
\end{array}  \]
Apply corollary~\ref{C:normmeaverage} for $(x_i)$, $\varepsilon$, $N=1$ and 
$M=j$ (and appropriate $(\varepsilon^k_i)$) to obtain $x$, a $(j,
(\varepsilon^k_i),1)$ average of $(x_i)$ w.r.t. $(e_i)$ with admissible 
averaging tree $(x^k_i)_{k=0,i=1}^{j,N^k}$ of $(x_i)$ and maximum coordinates
$(N^k_i)_{k=1,i=0}^{j,N^k}$ w.r.t. $(e_i)$. For $i=1, \ldots , N^0$ if $x^0_i=
x_s$ then define $j_i=\bar{j}_s$. Then $j_1< \cdots <j_{N^0}$ and for $i=1, 
\ldots , N^0$ we have
\[ \begin{array}{ll} 
 \forall p=1, \ldots , j_i, & \| x^0_i \|_{N^0_{i-1},p}< 
        \phi_p(1 + \varepsilon) \mbox{ and } \\
 \forall p \geq j_{i+1}, & \| x^0_i \|_{N^0_{i-1},p}< 
        \varepsilon. 
\end{array}  \]
Note (by remark~\ref{R:meaverage} (2)) 
that $x^j_1$ is $j$-admissible w.r.t. $(x_i)$ and by 
corollary~\ref{C:normmeaverage} (2) there exist $m \in \N$,
intervals $F_1< \ldots < F_m$ which don't split the $x^0_s$'s, and  $(p_\subl)
_{\ell=1}^m \subset \N$ with $p_\subl \geq j$ for all $\ell$ such that 
\[ \| x \| \leq \max \left( \bigcup_{\ell=1}^m \{ \frac{\theta
_{p_\subl}}{\theta_{p_\subl-j}} \| x^0_s \|_{N^0_{s-1},p_\subl-j}: 
\ran(x^0_s) \subseteq F_\subl \}^* \right) + \varepsilon. \]
For each $\ell=1, \ldots , m$ if $p_\subl >j$ then there exists 
exactly one $m_\subl \in \N$
such that $j_{m_\subl} \leq p_\subl-j < j_{m_\subl+1}$. We shall use the obvious
remark  that if $A \subseteq [0, \infty)$ is a finite non-empty set and $a \in 
A$ then $\max(A^*) \leq \max(A \backslash \{ a \})$. If $p_\subl =j$ then 
$\theta_{p_\subl}/ \theta_{p_\subl -j} =\theta_j$ and we note that 
$\| x^0_s \|_{N^0_{s-1},0} \leq \frac{1}{\theta_1} \| x^0_s \| = \frac{1}{
\theta_1}$. Thus
\[ \| x \| \leq \max \left( \bigcup_{ \{ \ell : p _\subl > j\} }
 \{ \frac{\theta_{p_\subl}}{\theta_{p_\subl-j}} \| x^0_s
 \|_{N^0_{s-1}, p_\subl-j}: 
\ran(x^0_s) \subseteq F_\subl, s \not = m_\subl \} \cup \{ \frac{\theta_j}{
\theta_1} \} \right) + \varepsilon. \]
Let $\ran(x^0_s) \subseteq F_\subl$ and $p_\subl >j$. If $s<m_\subl$ we have 
$j_{s+1} \leq j_{m_\subl} \leq p_\subl -j$ and so
$\| x^0_s \|_{N^0_{s-1},p_\subl-j}< \varepsilon$. If $s>m_\subl$ we have 
$j_s \geq j_{m_\subl +1}>p_\subl -j$ and so $\| x^0_s \|_{N^0_{s-1},p_\subl-j}<
 \phi_{p_\subl-j} ( 1+ \varepsilon)$. Note that 
\[ \frac{\theta_{p_\subl}}{\theta_{p_\subl-j}} \phi_{p_\subl-j}= \theta^j 
\phi_{p_\subl} \]
and therefore 
\[ \| x \| \leq \theta^j \sup_{p \geq j} \phi_p (1 + \varepsilon)\vee
\frac{\theta_j}{\theta_1} + 2 \varepsilon. \] 
Note (by remark~\ref{R:meaverage}) that we can write $x= \sum_F a_i x_i$ 
for some set $F \in S_j$ where $a_i > 0$ for all $i \in F$ and $\sum_{i \in F}
a_i=1$. Therefore $\delta_j(Y) \leq \| x \|$ and since $\varepsilon >0$ 
is arbitrary we obtain the result.  \hfill  $\Box$ 

\bigskip
\noindent
Note that theorem~\ref{T:computedeltaj} does not necessarily give the best 
possible estimate
for $\delta_j(Y)$. Indeed if $\theta_n=2^{-n}$ for all $n$ then $T=T(\theta_n,
S_n)_{\N}$ and for all $Y \prec T$, $\delta_j(Y)=2^{-j}$ \cite{OTW}. Yet
theorem~\ref{T:computedeltaj} only gives $\delta_j(Y) \leq 2^{-j+1}$. However
we have the following estimate which does yield the proper estimate for
Tsirelson's space.

\begin{Thm}
Let $X=T(\theta_n,S_n)_{\N}$ be regular. Then for all $Y \prec X$ and $j \in \N$
we have
\[ \delta_j(Y) \leq \theta^j \sup_{p \geq j} \frac{\phi_p}{\phi_{p-j}}.\]
\end{Thm}

\noindent
{\bf Proof} Let $Y \prec X$, $j \in \N$ and $\varepsilon >0$.
Since $Y$ contains $\ell_1^n$'s uniformly, for all $N \in \N$ $\exists y \in Y$
with $1= \| y \| \leq \| y \|_{N,0} \leq 1+ \varepsilon$. (see eg \cite{OTW} 
proposition 2.7). Therefore we may choose inductively a normalized block 
sequence $(x_i)$ in $Y$ so that for $i \in \N$ if $N_i = \max(\ran(x^0_i))$
w.r.t. $(e_i)$ ($N_0=1$) then $\| x_i \|_{N_{i-1},0} \leq 1 +
\varepsilon$. Note then that for every $i,p \in \N$, $\|x_i \|_{N_{i-1},p}
\leq \| x_i \|_{N_{i-1},0} \leq 1 + \varepsilon$. Apply 
corollary~\ref{C:normmeaverage} (for an appropriate sequence 
$(\varepsilon^k_i)$) to 
obtain  $x$ a $(j, (\varepsilon^k_i), 1)$ average of $(x_i)$ w.r.t. $(e_i)$ 
with admissible averaging tree $(x^k_i)_{k=0,i=1}^{j,N^k}$ of $(x_i)$ and 
maximum coordinates $(N^k_i)_{k=0,i=0}^{j,N^k}$ w.r.t. $(e_i)$. Note then that
for every $i,p \in \N$ we have that $\| x^0_i \|_{N^0_{i-1},p} \leq 1+ 
\varepsilon$. By corollary~\ref{C:normmeaverage} (2)
there exist $m \in \N$, $F_1< \cdots < F_m$ intervals in $\N$ which 
don't split the $x^0_i$'s and integers $(p_\subl)_{\ell=1}^m$ with $p_\subl 
\geq j$ for all $\ell$, such that 
\[ \| x \| \leq \max \{ \frac{\theta_{p_\subl}}{\theta_{p_\subl -j}}
\| x^0_i \|_{N^0_{i-1},p_\subl -j}: \ell =1, \ldots , m, \ran(x^0_i) \subseteq 
F_\subl \} + \varepsilon \leq \theta^j \sup_{p \geq j }\frac{\phi_p}{\phi_{p-
j}} (1+ \varepsilon) + \varepsilon \]
and the result follows since $\varepsilon >0$ is arbitrary. \hfill $\Box$

\bigskip
\noindent
To estimate $\delta_j(Y)$ for $Y=X$ is easy as we see from the next

\begin{Thm} 
Let $X=T(\theta_n,S_n)_{n \in \N}$ be regular. Then for all $j \in \N$ we have
$\delta_j(X)= \theta_j$.
\end{Thm}

\noindent
{\bf Proof} Let $j \in \N$ and $\varepsilon >0$. Apply 
corollary~\ref{C:normmeaverage} for $(x_i)=(e_i)$, $M=j$, $N=1$ and an
appropriate
sequence $(\varepsilon^k_i)$, to obtain $x$, a $(j,(\varepsilon^k_i),1)$ average
of $(e_i)$ w.r.t. $(e_i)$ with admissible averaging tree $(x^k_i)_{k=0, i=1}
^{j,N^k}$ and maximum coordinates $(N^k_i)_{k=0,i=0}^{j,N^k}$ w.r.t.
$(e_i)$. Then by (2) there exists $m \in \N$, $F_1< \cdots <F_m$ intervals in 
$\N$ and integers $(p_\subl)_{\ell =1}^m$ with $p_\subl \geq j$ for all $\ell$,
such that 
\[ \| x \| \leq \max \{ \frac{\theta_{p_\subl}}{\theta_{p_\subl -j}}
\| x^0_i \|_{N^0_{i-1},p_\subl -j}: \ell =1, \ldots , m, \ran(x^0_i) \subseteq 
F_\subl \} + \varepsilon. \]
Since $(x^0_i)_{i=1}^{N^0}$ is a subsequence of $(e_i)$, we have $\| x^0_i \|
_{N^0_{i-1},p_\subl -j} = \theta_{p_\subl -j}$ for every $i=1, \ldots , N^0$ 
and $\ell=1, \ldots , m$. Thus $\| x \| \leq \max_{1 \leq \ell 
\leq m} \theta_{p_\subl} + \varepsilon$. Since the sequence $(\theta_i)$ is
decreasing we have $\| x \| \leq \theta_j + \varepsilon$. Since $\supp(x) \in 
S_j$ and $\varepsilon >0$ is arbitrary we obtain the result. \hfill $\Box$

\bigskip
\noindent
{\bf Question} If $X= T(\theta_n,S_n)_\N$ is a regular mixed Tsirelson space
and $Y \prec X$ is $\delta_j(Y)= \theta_j$ for every $j \in \N$?

\noindent
G. Androulakis, Dept of Math, RLM 8.100,  The University of Texas at Austin, 
Austin, TX 78712 \\
e-mail: giorgis@math.utexas.edu 

\bigskip
\noindent
E. Odell, Dept of Math, RLM 8.100, The University of Texas at Austin, 
Austin, TX 78712 \\
e-mail: odell@math.utexas.edu


\begin{thebibliography}{99}

\bibitem[AA]{AA}
    D. Alspach, S. Argyros, {\em Complexity of weakly null sequences }, 
    Dissertationes Mathematicae {\bf 321} (1992).

\bibitem[AD]{AD} S. Argyros, I Deliyanni, {\em Examples of asymptotic $\ell_1$ 
    Banach spaces }, preprint.

\bibitem[FJ]{FJ} T. Fiegel, W. B. Johnson, {\em A uniformly convex Banach space
     which contains no $\ell_p$ }, Compositio Math. {\bf 29} (1974), 179-190.

\bibitem[Ma]{Ma} B. Maurey, {\em A remark about distortion}, Oper. Theory:Adv.
    Appl. {\bf 77} (1995), 131-142.

\bibitem[MMT]{MMT} B. Maurey, V. D. Milman, N. Tomczak-Jaegermann, {\em 
    Asymptotic infinite dimensional theory of Banach spaces}, Operator Theory:
   Advances and Applications, {\bf 77} (1995), 149-175.

\bibitem[MT]{MT} V. D. Milman, N. Tomczak-Jaegermann, {\em Asymptotic $\ell_p$
       spaces and bounded distortions}, Contemp. Math. {\bf 144} (1993), 
      173-196.

\bibitem[OS1]{OS1} E. Odell, Th. Schlumprecht, {\em The distortion problem},
       Acta Math. {\bf 173} (1994), 259-281.

\bibitem[OS2]{OS2} E. Odell, Th. Schlumprecht, {\em The distortion problem},
       GAFA {\bf 3} (1993), 201-207.

\bibitem[OS3]{OS3} E. Odell, Th. Schlumprecht, {\em Distortion and stabilized 
      structure in Banach spaces; New geometric phenomena for Banach and Hilbert
      spaces}, Proc. Inter. Cong. Math., Birkh\"{a}user Verlag, Basel (1995), 
      955-965.

\bibitem[OTW]{OTW} E. Odell, N. Tomczak-Jaegermann, R. Wagner {\em Proximity to
    $\ell_1$ and Distortion in asymptotic $\ell_1$ spaces}, preprint.

\bibitem[S]{S} Th. Schlumprecht, {\em An arbitrarily distortable Banach space},
     Israel J. Math. {\bf 76} (1991), 81-95.

\bibitem[T]{T} N. Tomczak-Jaegermann, {Banach spaces of type $p$ have 
     arbitrarily distortable subspaces}, preprint.

\bibitem[Ts]{Ts} B. S. Tsirelson, {\em Not every Banach space contains 
     $\ell_p$ or $c_{0}$}, Funct. Anal. Appl. {\bf 8} (1974), 138-141.
  
\end{thebibliography}
\end{document}